\documentclass[oneside,11pt,reqno]{amsart}
\usepackage[utf8]{inputenc}

\usepackage{amsfonts}
\usepackage{amsmath}
\usepackage{amssymb}
\usepackage{amsthm}
\usepackage[a4paper]{geometry}
\usepackage{mathrsfs} 
\usepackage[dvipsnames]{xcolor}

\usepackage{bbm}
\usepackage{mathtools}
\usepackage{dsfont}

\usepackage[pagebackref=true,colorlinks=true, linkcolor=RoyalPurple, citecolor=RoyalPurple]{hyperref}
\renewcommand*\backref[1]{\ifx#1\relax \else (Cited on #1) \fi}

\usepackage{appendix}
\usepackage{enumitem}
\usepackage{float}

\usepackage{tikz}
\usepackage{tikz-cd} 
\usetikzlibrary{arrows, arrows.meta}

\usepackage[nameinlink]{cleveref} 
\theoremstyle{plain}
\newtheorem{definition}{Definition}
\newtheorem{proposition}[definition]{Proposition}
\newtheorem{lemma}[definition]{Lemma}
\newtheorem{corollary}[definition]{Corollary}
\newtheorem{theorem}[definition]{Theorem}

\newtheorem{remark}[definition]{Remark}

\theoremstyle{definition}

\numberwithin{definition}{section}
\numberwithin{equation}{section} 

\newcommand*{\R}{\mathbb{R}}

\newcommand*{\N}{\mathbb{N}}

\newcommand*{\Fcal}{\mathcal{F}}

\newcommand*{\Pcal}{\mathcal{P}}
\newcommand*{\Gcal}{\mathcal{G}}

\newcommand*{\rhoc}{\rho_{\rm c}}

\newcommand{\abs}[1]{\left\lvert #1 \right\rvert}

\newcommand*{\1}{\mathds{1}}
\newcommand{\x}{\boldsymbol{x}}
\newcommand*{\y}{\boldsymbol{y}}

\newcommand{\oo}{\boldsymbol{o}}
\renewcommand*{\L}{\Lambda}

\renewcommand*{\d}{\mathrm{d}}
\newcommand*{\e}{\mathrm{e}}

\newcommand*{\ar}{{\rm ar}}
\newcommand*{\hc}{{\rm hc}}

\crefname{equation}{}{}

\title[Infinite-range multibody Gibbs variational principle]{The variational principle for a marked Gibbs point process with infinite-range multibody interactions}
\author[B. Jahnel]{Benedikt Jahnel}
\author[J. K\"oppl]{Jonas K\"oppl}
\author[Y. Steenbeck]{Yannic Steenbeck}
\author[A. Zass]{Alexander Zass}

\address[Benedikt Jahnel]{TU Braunschweig, Institut für Mathematische Stochastik, 
Germany, and Weierstrass Institute.}
\email{jahnel@wias-berlin.de}
\address[Jonas K\"oppl]{Weierstrass Institute, 
Germany.}
\email{koeppl@wias-berlin.de}
\address[Yannic Steenbeck]{TU Braunschweig, Institut für Mathematische Stochastik, 
Germany.}
\email{yannic.steenbeck@tu-braunschweig.de}
\address[Alexander Zass]{Weierstrass Institute, 
Germany.}
\email{zass@wias-berlin.de}
\date{\today}


\definecolor{amethyst}{rgb}{0.6, 0.4, 0.8}

\usepackage[pagewise]{lineno}

\begin{document}

\begin{abstract}
We prove the Gibbs variational principle for the Asakura--Oosawa model in which particles of random size obey a hardcore constraint of non-overlap and are additionally subject to a temperature-dependent area interaction. The particle size is unbounded, leading to infinite-range interactions, and the potential cannot be written as a $k$-body interaction for fixed $k$.
As a byproduct, we also prove the existence of infinite-volume Gibbs point processes satisfying the DLR equations. The essential control over the influence of boundary conditions can be established using the geometry of the model and the hard-core constraint.
\end{abstract}

\maketitle

\section{Introduction}
The Gibbs variational principle is a cornerstone in the equilibrium theory of spatially interacting particle systems.
In a nutshell, according to the laws of thermodynamics, equilibrium states of a system (subject to thermodynamic principles) must minimize the excess free energy, which reflects the competition between minimizing energy and maximizing entropy. 

On the other hand, equilibrium systems of interacting particles are often described by their conditional finite-volume distributions, in terms of Boltzmann weights with respect to some boundary conditions. Under a consistency condition on these finite-volume measures, one can then hope to find infinite-volume measures that are still consistent with the finite-volume description, and are able to properly describe the macroscopic behavior of the system. This infinite-volume consistency is described via the so-called DLR equations, and ensuring the existence of such measures is one the key tasks in statistical mechanics. In particular, the DLR equations can potentially be solved by multiple infinite-volume measures; this non-uniqueness of equilibrium states constitutes a type of first-order phase transition. However, it is not automatic for a measure to solve the DLR equations and also minimize the excess free energy.

In the case of lattice systems, under the assumption of translation invariance, this link between the DLR equations and the thermodynamic description of equilibria has been established under quite general assumptions for a very large class of models, see~\cite{georgii_gibbs_2011,preston1976random,lanford1969observables}. However, moving towards continuum systems, where particles are placed in Euclidean space and even particle numbers in fixed volumes are allowed to fluctuate, the literature becomes much more sparse.
In~\cite{Georgii_1994,georgii_equivalence_1995} the Gibbs variational principle is derived for bounded and unbounded pair potentials in the context of large-deviation results, and with the goal to establish equivalence of ensembles. In particular, no higher-order interactions are allowed, and the particles may not carry individual marks that would for example allow to model interactions between particles of varying sizes.
Next, in~\cite{dereudre2009variational} the Gibbs variational principle is established for a particular class of planar models in which (unmarked) particles interact via their associated Delaunay triangulation.
More recently,~\cite{dereudre_variational_2016} considers the case of Gibbs point processes with a finite-range interaction and where the existence of the energy density is guaranteed. However, establishing the existence of the energy density is not always easy and needs to be shown specifically for the model at hand. Moreover, many interesting models do not have a finite-range interaction, see, for example,~\cite{Dereudre_2009,dereudre_vasseur_2020,Dereudre_2023,Kondratiev2012,Lewin_2022,Roelly_Zass_2020,Zass_2022} and the references therein.

The contribution of the present manuscript is therefore to provide a first example of a Gibbs variational principle in the continuum for a system of particles with unbounded marks and with a multibody infinite-range interaction in which also an unbounded number of particles may jointly interact. 

We focus on the Asakura--Oosawa model with depletion interaction~\cite{AO54,asakura1958interaction,vrij1976polymers}, which may be informally described as follows. Imagine a volume with a random number of balls of a random size. The balls cannot penetrate each other, leading to a hardcore repulsion. Additionally, there are second-class particles with a fixed size which may penetrate each other without restriction, but not the balls of the first type. Integrating out the Poisson point process of second-class particles, it can be observed that the balls of type one must additionally obey an area interaction: for them to take up more space comes with the additional cost that, in that space, no second-class particles are present. For more details on the relevance of this 
colloid-polymer mixture model, see~\cite{vrij1976polymers,binder2014perspective,roth2000depletion,lekkerkerker2024colloids}. Let us mention that this model can be seen as a generalization of the Widom--Rowlinson model~\cite{widom1970new,chayes1995analysis,dereudre_houdebert_2019,dereudre2021sharp,jahnel2017widom} with an additional internal hardcore. The existence of infinite-volume Gibbs measures (in the sense of the DLR equations) for the Asakura--Oosawa model can be deduced from the criteria presented in~\cite{Roelly_Zass_2020}. Recently, the model has been analyzed more extensively, for example via cluster-expansion techniques~\cite{jansen2020cluster,jansen2020mayer,wittmann2022geometric} and with respect to the corresponding Langevin dynamics~\cite{FR23,fradon2025}. 

The proof rests on establishing existence of energy and entropy densities for which it is essential to control the influence of far away particles. This is usually achieved by a finite-range or deterministic decay assumption on the (pair) interaction potential. However, in our case we work only with an assumption on the tails of the mark distribution, and use the hard-core constraint as well as geometric arguments.  

The manuscript is organized as follows. 
In Section~\ref{sec:setting-and-notation}, we introduce the framework in which we are working, and set up the required notation before we state our main results in Section~\ref{sec_results}. We then proceed by explaining the structure and key ideas of the main proof in Section~\ref{sec_proof_strategy}. After this, we finally start with the main work and carry out the proof of our main results in Sections~\ref{sec_preliminaries}--\ref{sec_variational_principle}. 
    
\section{Setting and Notation}\label{sec:setting-and-notation}
 Throughout the article, we will often use the short-hand notation $f\lesssim g$ to express that there exists a finite constant $c>0$, independent of $f$ and $g$, such that $f\le c \cdot g$. 
 We write $f\lesssim_A g$ if $c=c(A)$ depends on some parameter $A$, for example the dimension $d$ of the ambient space.
\subsection{Point process formalism}
    Consider the space $\Omega$ of locally finite measures on the state space $\R^d\times [0, \infty)$, $d\geq 1$, that is 
    $$
        \Omega:=\Big\{\omega= \sum_{i\ge 0}\delta_{\x_i} \ \vert\ \x_i=(x_i,R_i)\in\R^d\times[0, \infty)\Big\}.
    $$
    Since, in what follows, we consider only simple point configurations, we may identify them with their support, i.e., $\omega = \sum_{i\ge 0} \delta_{\x_i} \equiv \{\x_i\}_{i\ge 0}\subset\R^d\times[0, \infty)$. For any subset $\L\subseteq\R^d$, we denote by $\omega_\L$ the restriction of the configuration $\omega\in\Omega$ to the set $\L\times[0, \infty)$, that is $\omega_\L = \omega\cap(\L\times[0, \infty))$. 
    Furthermore, for Borel sets \(\L \in \mathfrak{B}(\mathbb{R}^d)\), \(A \in \mathfrak{B}([0,\infty))\), and positive measurable functions \(f \colon \mathbb{R}^d \times [0, \infty) \to [0, \infty]\), we define the counting variables   
    \begin{align*}
        N_{\L, A}^{f}(\omega)
        := \sum_{\x = (x, R_x) \in \omega \cap (\L \times A)} f(x, R_x),
    \end{align*} 
    and denote in particular
    \begin{align*}
        N_{\L}^{f} := N_{\L, [0, \infty)}^{f}, \quad
        N_{\L, A} := N_{\L, A}^{1}, \quad\text{ and }\quad
        N_{\L} := N_{\L}^{1}.
    \end{align*}
    Next, we endow $\Omega$ with the canonical $\sigma$-algebra \(\mathcal{F}\) generated by the family of local counting variables, that is
    \begin{align*}
        \mathcal{F} = \sigma\big(N_{\L, A} \colon \L \in \mathfrak{B}_b(\mathbb{R}^d), A \in \mathfrak{B}_b([0, \infty)) \big),
    \end{align*} where \(\mathfrak{B}_b\) denotes the set of bounded Borel sets.
    We also fix some special sets of the position space \(\mathbb{R}^d\), which will be denoted by 
    \begin{align*}
        \Lambda_n := \left[-n/2, n/2\right)^d, \ n \in \mathbb{N},\quad\text{ and }\quad 
        C := \Lambda_1
    \end{align*} 
    throughout the whole manuscript.
    
    In order to define classes of test functions, let \(\psi \colon [0, \infty) \to [1, \infty)\) be given by
    \begin{align*}
        \psi(R) := 1 + R^d, \quad R \in [0, \infty),
    \end{align*}
    and call a function \(f \colon \Omega \to \mathbb{R}\)  {\em local}, if \(f = f(\cdot_{\L})\)  and {\em tame} if \(\abs{f} \lesssim (1 + N_\L^{\psi})\) for some \(\L \in \mathfrak{B}_b(\mathbb{R}^d)\).  The set of all tame and local functions will be denoted by \(\mathcal{L}\).
    
    Finally, we denote by $\Pcal$ the set of all probability measures \(P\) on $\Omega$, and by $\Pcal_\Theta$ the set of all translation-invariant probability measures on $\Omega$ with finite \(\psi\)-intensity, i.e., with \(P[N^{\psi}_C] < \infty\). Here and in the sequel, $P[f]$ denotes the expected value of $f$ w.r.t.~the probability measure $P$.
    We equip \(\Pcal_\Theta\) with the so-called \(\tau_\mathcal{L}\)-topology, which is the coarsest topology that makes the mappings \(P \mapsto P[f]\) continuous for all \(f \in \mathcal{L}\).

\subsection{Colloidal multi-body interaction}
    The following model describes the effective interaction between colloids in a polymer bath; it is known as {\em depletion interaction} in the Asakura--Oosawa model, see e.g.~\cite{AO54}.
    For a finite configuration $\omega = \{(x_i,R_i)\}_{1\leq i\leq n} = \{\x_i\}_{1\leq i\leq n}$, the interaction is given by
    \begin{equation}\label{eq_1}
        H(\omega) := H^\hc(\omega) + H^\ar(\omega), 
    \end{equation} 
    where
    \begin{itemize}
        \item $H^\hc(\omega) := \sum_{1\leq i<j\leq n} \phi_\hc(\x_i,\x_j) = (+\infty)\sum_{1\leq i<j\leq n} \1\{\abs{x_i-x_j}\leq R_i+R_j\}$ is a hard-core pair interaction; 
        \item $H^\ar(\omega) := \big\lvert \bigcup_{i=1}^n B(x_i,R_i+r) \big\vert$ is the total area covered by $\omega$, with radii enlarged by $r>0$. 
    \end{itemize}
    In other words, the interaction models colloidal particles of radius $R_i$ centered at $x_i$ in a bath of polymers of size $r$.
    Note that, by inclusion-exclusion, the area-interaction term can be rewritten in terms of $k$-body potentials as
    \begin{equation*}
        H^\ar(\omega) = \sum_{k\geq 1} \,\, \sum_{1\leq i_1 < \ldots < i_k \leq n} \phi_k(\x_1,\ldots,\x_k),
    \end{equation*}
    where
    \begin{equation*}
     \phi_k(\x_1,\dots,\x_k) := (-1)^{k-1}\big\lvert B(x_1,R_1+r)\cap\ldots\cap B(x_k,R_k+r)\big\rvert.
    \end{equation*}

    \begin{remark}[$k$-body interactions via sufficiently small $r$]
        Let us note that there exists a value $\rhoc:=\rhoc(d)>0$ such that, if $r/R\leq\rhoc$, where $R=\min\{R_1,\dots, R_n\}$, the $k$-body interactions vanish for $k\geq 4$, that is,
        \begin{equation}\label{eq:3body}
            H^\ar(\omega)=\sum_{i=1}^n \phi_1(\x_i) +\sum_{1\leq i<j\leq n}\phi_2(\x_i,\x_j) + \sum_{1\leq i<j<k\leq n} \phi_3(\x_i,\x_j,\x_k).
        \end{equation}
        In dimension $d=2$, this value is obtained by considering four spheres, with maximal radius $R$ centered on the vertices of a square of side length $2R$. There can be a $4$-way overlap only if the distance between any vertex and the center of the square is less than $2(R+r)$, that is if $r/R>\rhoc(2)= \sqrt{2}-1$. An analogous reasoning applies in dimension $d\geq 3$. Replacing the square by a tetrahedron, this for example gives $\rhoc(3) = \sqrt{3/2}-1$. 
    \end{remark}

    \medskip
    Next, we define the {\em conditional energy} of $\omega \in \Omega$ in a bounded subset $\L\subseteq \R^d$ {\em with boundary configuration} \(\zeta \in \Omega\) via
    \begin{equation}
            H_{\Lambda, \zeta}(\omega) := H^\hc_{\Lambda, \zeta}(\omega)+H^\ar_{\Lambda, \zeta}(\omega)
        \end{equation}
        with $H^\hc_{\Lambda, \zeta}(\omega):= \sum_{\{\x,\y\}\subseteq \omega_\L \zeta_{\L^c} \colon \{\x,\y\}\cap\L\neq\emptyset} \phi_\hc(\x,\y)$ and 
    \begin{equation}
            H^{\rm ar}_{\Lambda, \zeta}(\omega) :=\Big\vert \bigcup_{\x \in \omega_{\Lambda}} B(x, R_x + r) \,\setminus\, \bigcup_{\y \in \zeta_{\Lambda^c}} B(y, R_y + r) \Big\vert,
    \end{equation} 
    where $\Lambda^c:=\R^d\setminus\Lambda$ and $\omega_\L \zeta_{\L^c} = \omega_\L\cup\zeta_{\L^c}$ denotes the concatenation of the configurations $\omega_\L$ and $\zeta_{\L^c}$. 
    For the choice \(\zeta = \emptyset\), we call \(H_{\L, \emptyset}(\omega)\) the ({\em conditional}) {\em energy} of $\omega$ in \(\L\) {\em with free boundary conditions}.
    Another noteworthy special case is if \(\zeta = \omega\) and we then just write \(H_{\L}(\omega)\) instead of \(H_{\L, \omega}(\omega)\) and speak of the {\em conditional energy} of \(\omega\) {\em in} \(\L\).

\subsection{Gibbsian formalism} 
    Based on the energy, we now define the finite-volume equilibrium Gibbs measures with respect to suitable boundary conditions. For this, let $\pi^z_\Lambda$ denote the iid-marked Poisson point process on $\Lambda \times [0,\infty)$ with intensity measure $z \1\{x\in \L\}\mathrm{d}x  \otimes \mathfrak{R}(\d R)$, with $\Lambda\in \mathfrak{B}(\R^d)$ and $\mathfrak{R}$ a probability measure on $[0,\infty)$. We write $\pi^z:=\pi^z_{\R^d}$. 
    Throughout the paper, we make the following integrability assumption, 
    $$\int \e^{R^{d+\delta}}\mathfrak{R}(\d R)<\infty\qquad \text{ for some }\delta>0.$$

    \begin{definition}[Finite-volume Gibbs measures]
        For any bounded $\Lambda\in \mathfrak{B}(\R^d)$, the {\em finite-volume Gibbs measure} with parameters \(z, \beta > 0\) and boundary condition $\zeta \in \Omega$ is the probability measure $G_{\Lambda, z, \beta, \zeta}$ on point configurations in $\Lambda$ given by
        \begin{equation}
            G_{\Lambda, z, \beta, \zeta}(\d\omega_\L) = \frac{1}{Z_{\L, z, \beta, \zeta}}\e^{- \beta H_{\Lambda, \zeta}(\omega_\L)}\pi^z_\Lambda(\d\omega_\L),
        \end{equation}
    \end{definition}
    
    \noindent 
    Let us note that, since $H_{\Lambda, \zeta} \geq 0$ and \(H_{\Lambda, \zeta}(\emptyset) = 0\), we have that $\exp(-z|\Lambda|)\le Z_{\L, z, \beta, \zeta} \le 1$ and thus the finite-volume Gibbs measures are well-defined for all boundary conditions. One can interpret the random radii $R_i$ as \textit{annealed} disorder in the Gibbs measure.

    \begin{definition}[Gibbs measures]\label{def_infinite_volume_Gibbs_measure}
        A probability measure $P\in\Pcal$ is called an {\em infinite-volume Gibbs measure} for \(z, \beta > 0\) if the following DLR equation holds for any bounded $\L\in \mathfrak{B}(\R^d)$ and any measurable function $f\geq 0$,        \begin{equation}\label{eq:DLR}
            \int P(\d\omega) \, f(\omega)
                	= \int P(\d\zeta) \int \, G_{\L, z, \beta, \zeta}(\d\omega_\L) \, f(\omega_{\Lambda}\zeta_{\Lambda^c}).
        \end{equation} 
        We write $\Gcal^{z,\beta}(H)$ for the set of all infinite-volume Gibbs measures for \(z, \beta > 0\). We denote by $\Gcal^{z,\beta}_\Theta(H)= \Gcal^{z,\beta}(H)\cap\Pcal_\Theta$ the set of all translation-invariant infinite-volume Gibbs measures.
    \end{definition}
    In the following, the infinite-volume Gibbs measure we consider will always be translation-invariant.
    
    Let us already note that our main result will in particular establish existence of infinite-volume translation-invariant Gibbs measures for the model defined in~\Cref{eq_1}. However, the existence could also be verified by checking the conditions presented in~\cite{Roelly_Zass_2020}.

    \subsection{Specific entropy}
    On the way to present our main result, the Gibbs variational principle, we need to introduce the three main ingredients of the variational formula, the {\em specific entropy}, the {\em energy density}, and the {\em pressure}. 
    
        First, we introduce the entropic part and recall its basic properties. For this, consider the (finite-volume) {\em relative entropy} for bounded $\L\in\mathfrak{B}(\R^d)$ given by
\begin{align*}
            I(P_{\L} \,\vert\, \pi^z_{\L})=\int P_{\L}(\d \omega_\L) \, \log\left(\frac{\d P_{\L}}{\d \pi_\L^z}(\omega_\L)\right),
        \end{align*}
    if $P_\L\ll\pi_\L^z$ and $I(P_{\L} \,\vert\, \pi^z_{\L})=\infty$ otherwise. The proof of the following standard result can be found for example in~\cite[Proposition 2.6]{GZ93}.
    \begin{proposition}[Specific entropy]\label[proposition]{prop_specific_entropy_basics}
        For every activity \(z > 0\) and \(P \in \mathcal{P}_\Theta\), the {\em specific entropy}
        \begin{equation}\label{eq_def_entropy_is_sup}
            I^z(P) := \lim_{n \uparrow \infty} \frac{I(P_{\L_n} \,\vert\, \pi^z_{\L_n})}{\abs{{\L_n}}}=\sup_{\L \in \mathfrak{B}_b(\mathbb{R}^d)} \frac{I(P_\L \,\vert\, \pi^z_\L)}{\abs{\L}}
        \end{equation}
        exists.
        Moreover, the map \(P \mapsto I^z(P)\) is lower-semicontinuous and the level sets \(\{I^z \leq c\}, c \in \mathbb{R},\) are compact and sequentially compact w.r.t.~the \(\tau_\mathcal{L}\)-topology.
    \end{proposition}

    \subsection{Energy density}
    The energetic component in the Gibbs variational formula, namely the \textit{energy density} is formally defined as 
    \begin{align*}
        H(P) := \lim_{n \uparrow \infty}\frac{P[H_{\Lambda_n \emptyset}]}{\abs{\Lambda_n}}. 
    \end{align*}
    In contrast to the specific entropy, its existence and continuity properties are not guaranteed by general principles and need to be established in a model-dependent way. This will already constitute one part of our main result. 

    \subsection{The pressure}
    The third ingredient of the variational principle is the so-called \textit{pressure}, i.e., the density limit of the partition function, formally defined as  
    \begin{align*}
        \psi(z,\beta) 
        :=
        \lim_{n \to \infty} \frac{\log(\pi^z_{\Lambda_n}\left[\e^{-\beta H_{\Lambda_n, \emptyset}} \right])}{\vert \Lambda_n \vert}.
    \end{align*}
    Here, the existence of the limit \(\lim_{n \to \infty} \frac{\log(\pi^z_{\Lambda_n}\left[\e^{-\beta H_{\Lambda_n, \zeta}} \right])}{\vert \Lambda_n \vert}\) is also quite subtle and depends on the choice of the boundary condition $\zeta$. We will show that, for a class of sufficiently nice boundary conditions (including the free case \(\zeta = \emptyset\)), it exists and actually does not depend on the boundary condition itself. 
    
\section{Main result}\label{sec_results}

The following Gibbs variational principle is the main result of this work.  
    \begin{theorem}[Gibbs variational principle]\label{thm_variational_principle}
        Let $z>0$ and $\beta > 0$ be fixed.  
        \begin{enumerate}
            \item[i.] The functional 
            \begin{align*}
                \mathcal{P}_\Theta \ni P \mapsto I^z(P) + \beta H(P) 
                \in [0,\infty]
            \end{align*}
            is well-defined and lower-semicontinuous.
            \item[ii.] The functional attains its finite infimum on $\mathcal{P}_\Theta$ and this infimum coincides with the pressure, i.e., 
            \begin{align*}
                \inf_{P \in \mathcal{P}_{\Theta}} \left[I^z(P) + \beta H(P)\right] 
                = 
                - \psi(\beta,z).
            \end{align*}
            \item[iii.] A measure $P \in \mathcal{P}_\Theta$ minimizes the functional if and only if it is a translation-invariant Gibbs measure, i.e., $P\in\Gcal^{z, \beta}_\Theta$. 
        \end{enumerate}
    \end{theorem}
Note again that, as a byproduct of our main result, existence is guaranteed.
\begin{corollary}[Existence]\label{prop_existence_of_infinite_volume_gibbs_measure}
        For any activity $z>0$ and any inverse temperature $\beta>0$, there exists at least one translation-invariant infinite-volume Gibbs measure for the model defined in~\Cref{eq_1}. 
\end{corollary}

    Before we shift gears and discuss the overall strategy of the proof, let us remark that  the statement of \Cref{thm_variational_principle} remains true if one replaces \(\mathcal{P}_\Theta\) with the set of all translation-invariant probability measures on \(\Omega\). The energy density \(H(P)\) is still well-defined and can even be infinite, but \(P[N_C] = \infty\) necessarily implies \(I^z(P) = \infty\). \medskip 

    Moreover, we conjecture that \Cref{thm_variational_principle} still holds if one makes minor modifications to the model such as:
    \begin{enumerate}
        \item replacing spheres by some other type of convex body,
        \item replacing the hard-core part by a Lennard--Jones-type soft-core potential (with sufficiently strong short-range repulsion), 
        \item using an interaction potential which is not the area interaction but satisfies the properties in \Cref{lem_energy} below.  
    \end{enumerate}  
    All of these generalizations seem reasonable, but one would certainly pay the price of longer and harder-to-read proofs and potentially obscuring the main ideas. Let us also mention that our integrability assumption on the marks is possibly not optimal and the Gibbs variational principle still holds under weaker assumptions.
\section{Strategy of proof}\label{sec_proof_strategy}
For the proof of \Cref{thm_variational_principle} we will treat the main ingredients of the variational principle, namely the energy density and the pressure, separately, and then conclude by putting everything together. 
Let us briefly explain some of the key steps and note that this can also serve as a blueprint for establishing Gibbs variational principles for other models.  
\subsection{The energy density and its continuity properties}
Our first result deals with the existence and continuity properties of the energy density. 
\begin{proposition}\label[proposition]{prop_existence_of_energy_density}
    For all $P \in \mathcal{P}_\Theta$ the energy density 
    \begin{align*}
        H(P) := \lim_{n \uparrow \infty} \frac{P[H_{\L_n, \emptyset}]}{\abs{\L_n}}
    \end{align*}
    exists. The functional \begin{align*}
        \mathcal{P}_\Theta \ni P \mapsto H(P) \in [0,\infty]
    \end{align*}
    is continuous on $\{P\in \mathcal P_\Theta\colon H(P) < \infty\}$ and in particular lower semicontinuous on $\mathcal{P}_\Theta$ with respect to the $\tau_{\mathcal{L}}$-topology. 
\end{proposition}

The proof is based on deriving an explicit representation of $H(P)$ in terms of the Palm measure of $P$. 
Proving the continuity of $H(\cdot)$ on $\{P\in \mathcal P_\Theta\colon H(P) < \infty\}$ requires some geometric considerations to control the unbounded-range interactions. The details are provided in \Cref{sec_energy_density}. 
Let us also note that the map $P \mapsto H(P)$ is not continuous on the whole set $\mathcal{P}_\Theta$. An example of point of discontinuity is given in \Cref{remark_discontinuity}. 

\subsection{Tempered configurations and measures}
Especially for the existence of the pressure, we need to control the effect of considering different boundary conditions. There, it will be very useful to only work with boundary conditions which are sufficiently nice, these are the so-called \textit{tempered} configurations. 
Of course, the set of nice configurations also needs to be sufficiently large such that all relevant measures are concentrated on it, so there is a certain trade-off. 
Let us start by introducing some more notation. Fix some \(\gamma \in (0, \delta/d)\) and set 
        \begin{align*}
            \epsilon := \Bigl(1 - \frac{\gamma d}{\delta}\Bigr)\frac{\delta}{d+\delta},
        \end{align*}
        where \(\delta > 0\) is such that \(\int_{[0, \infty)} \e^{r^{d+\delta}}\mathfrak{R}(\d r) < \infty\) and define
        \begin{align*}
            g(n) := n^{1 - \epsilon}, \quad n \in \mathbb{N}.
        \end{align*}
In particular, 
        \begin{align*}
            \mathfrak{R}(R \geq g(n))
            \lesssim_{\mathfrak{R}, \delta} \e^{-\abs{\L_n}^{1+\gamma}}.
        \end{align*}
Then, we define the set of {\em tempered configurations} by
        \begin{align*}
            \Omega^* = \left\{\omega \in \Omega \ \middle\vert \ \exists N \in \mathbb{N}\ \forall n \geq N \ \forall \x \in \omega_{\Lambda_n} \colon R_x \leq g(n)  \right\}.
        \end{align*}
With this notation at hand, we can now introduce the class of measures we will be particularly interested in. 
\begin{definition}[Temperedness]
        A measure $P \in \mathcal{P}_\Theta$ is called {\em tempered}, if 
        \begin{align*}
            P(\Omega^*) = 1.
        \end{align*}
\end{definition}

Restricting our attention to tempered measures is not problematic for the following two reasons. First, non-tempered measures have infinite specific entropy, so they will not appear as minimizers in the Gibbs variational principle.   

\begin{lemma}[Finite specific entropy implies temperedness]\label{lemma_support_of_measure_with_finite_free_energy}
        Let \(P \in \mathcal{P}_{\Theta}\) with \(I^z(P) < \infty\), then \(P\) is tempered.
    \end{lemma}

Second, all infinite-volume Gibbs measures for the model we consider are tempered.

\begin{lemma}[Infinite-volume Gibbs measures are tempered]\label{le_infinite_volume_gibbs_measures_are_tempered}
        Let \(P \in \mathcal{P}_{\Theta}\) be an infinite-volume translation-invariant Gibbs measure \(P \in \Gcal_\Theta^{z, \beta}(H)\), then \(P\) is tempered.
\end{lemma}

The proofs of \Cref{lemma_support_of_measure_with_finite_free_energy} and \Cref{le_infinite_volume_gibbs_measures_are_tempered} can be found in \Cref{sec_temperedness}. 
Summarizing, one sees that only considering tempered measures is not really a restriction but enables us to proceed as follows. 

\subsection{Existence of the pressure}
The next step is to show the existence of the pressure, i.e., the density limit of the partition function. We do this step-by-step by exchanging the boundary conditions and start out by considering the periodic case. For this, we define $H_{\L_n, \mathrm{per}}(\omega_{\L_n})$ to be the finite-volume energy when $\L_n$ is identified with the $d$-dimensional torus. 

\begin{proposition}\label[proposition]{prop_pressure_per}
            We have
            \begin{align*}
                \lim_{n \to \infty} \frac{\log\big(\pi^z_{\Lambda_n}\big[\e^{-\beta H_{\Lambda_n, \mathrm{per}}} \big]\big)}{\vert \Lambda_n \vert}
                = -\inf_{P \in \mathcal{P}_{\Theta}} \left[I^z(P) + \beta H(P)\right].
            \end{align*}
        \end{proposition}

        Let us now write $H_{\L,\emptyset}(\omega) :=H_{\L}(\omega_\L)$ for the energy with free boundary conditions.
        Extending \Cref{prop_pressure_per} to free boundary conditions is essentially a consequence of the following comparison between free and periodic boundary conditions.
    \begin{proposition}\label[proposition]{le_comparison_pressure_per_and_free}
        We have
        \begin{align*}
            \lim_{n \to \infty} \frac{\log\big(\pi^z_{\L_n}\big[\e^{-\beta H_{\L_n, \mathrm{per}}}\big]\big)}{\vert \L_n \vert}
            = \lim_{n \to \infty} \frac{\log\big(\pi^z_{\L_n}\big[\e^{-\beta H_{\L_n, \emptyset}} \big]\big)}{\vert \L_n \vert}.
        \end{align*}
    \end{proposition}

    After this intermediate step, one can then extend the existence of the pressure to boundary conditions sampled from a tempered measure $P$ with finite energy density. 
    
    \begin{proposition}\label[proposition]{prop_pressure_with_boundary_conditions}
            Let \(P \in \mathcal{P}_{\Theta}\) be tempered and with \(H(P) < \infty\). Then,
            \begin{align*}
                \lim_{n \to \infty} \frac{\log\big(\pi^z_{\Lambda_n}\big[\e^{-\beta H_{\Lambda_n, \zeta}} \big]\big)}{\vert \Lambda_n \vert}
                = -\inf_{P \in \mathcal{P}_{\Theta}} \left[I^z(P) + \beta H(P)\right]
            \end{align*} for \(P\)-a.e.~\(\zeta \in \Omega\).
        \end{proposition}
    The proof is again based on a comparison estimate, now with the case of free boundary conditions. All the details can be found in \Cref{sec_pressure}. 
\subsection{The variational principle}
In the end, we put all of the previous ingredients together to prove our main result \Cref{thm_variational_principle}. For this, our strategy is based on the two-body situation in \cite{georgii_equivalence_1995}. 
As a first step, we show that the variational principle is non-trivial in the following sense. 
\begin{lemma}\label[lemma]{le_free_energy_attains_finite_infimum}
        We have 
        \begin{align*}
            \inf_{P \in \mathcal{P}_{\Theta}} \left[I^z(P) + \beta H(P)\right]
            = \min_{P \in \mathcal{P}_{\Theta}} \left[I^z(P) + \beta H(P)\right]
            < \infty,
        \end{align*} 
        i.e., \(P \mapsto I^z(P) + \beta H(P)\) attains its infimum and it is finite.
    \end{lemma}

    To show that the set of minimizers coincides with the translation-invariants infinite-volume Gibbs measures, we need to compare measures $P \in \mathcal{P}_\Theta$ with finite-volume Gibbs measures. For this note that in finite volumes $\Lambda_n$ we have
    \begin{align*}
        \log(Z_{\L_n, z, \beta}) + \beta P_{\L_n}[H_{\L_n, \emptyset}] +  I(P_{\L_n} \,\vert\, \pi^z_{\L_n}) =  P_{\L_n}\left[\log\left(\frac{\mathrm{d}P_{\L_n}}{\mathrm{d}G_{\L_n, z, \beta}}\right)\right]
        =: I(P_{\L_n} \,\vert\, G_{\L_n, z, \beta}),
    \end{align*}
    where $G_{\Lambda_n,z,\beta}$ is just the finite-volume Gibbs measure with free boundary conditions and normalization constant $Z_{\Lambda_n, z, \beta}$. 
    Of course, the term on the right-hand side is nothing but a relative entropy.
    By using that the energy density and the pressure exist one sees that the relative entropy density exists and satisfies 
    \begin{align*}
        \delta F(P) := \lim_{n \to \infty}\frac{I(P_{\L_n} \,\vert\, G_{\L_n, z, \beta})}{\abs{\Lambda_n}}
        =
        \left[I^z(P) +\beta H(P)\right]-\inf_{Q \in \mathcal{P}_\Theta}\left[I^z(Q)+\beta H(Q)\right]. 
    \end{align*}
    So a measure $P \in \mathcal{P}_\Theta$ is a minimizer, if and only if its relative entropy density (or excess free energy) vanishes. By controlling the  effect of different boundary conditions this identity allows us to show that every infinite-volume translation-invariant Gibbs measure is indeed a minimizer. 

    \begin{proposition}\label[proposition]{thm_gibbs_are_minimizers}
            Let \(P\) be a translation-invariant infinite-volume Gibbs measure for activity $z$ and inverse temperature $\beta$, i.e., \(P \in \Gcal_\Theta^{z, \beta}(H)\),
            then 
            \begin{align*}
                I^z(P) + \beta H(P)
                =
                \inf_{Q \in \mathcal{P}_{\Theta}} \left[I^z(Q) + \beta H(Q)\right].
            \end{align*}
    \end{proposition}  

    As a final step, we need to argue that the translation-invariant Gibbs measures are actually the only minimizers. As usual, this is the harder part of the Gibbs variational principle. 
    
    \begin{proposition}\label[proposition]{thm_minimizers_are_gibbs}
            Let \(P \in \mathcal{P}_\Theta\).
            If 
            \begin{align*}
                I^z(P) + \beta H(P)
                =
                \inf_{Q \in \mathcal{P}_{\Theta}} \left[I^z(Q) + \beta H(Q)\right],
            \end{align*}
            then \(P\) is a translation-invariant infinite-volume Gibbs measure for activity $z$ and inverse temperature $\beta$, i.e., \(P \in \Gcal_\Theta^{z, \beta}(H)\).
        \end{proposition}
        
    Recall that, in order to show that $P \in \mathcal{P}_\Theta$ is a Gibbs measure, one needs to show that for any bounded measurable $\L \subset \R^d$ and any measurable $f\geq 0$, one has 
     \begin{align*}
        \int \, P(\mathrm{d}\zeta) \, f(\zeta) 
            = 
        \int \, P(\mathrm{d}\zeta) \, \left[\int \, G_{\L, z, \beta, \zeta}(\mathrm{d}\omega) \, f(\omega_{\L}\zeta_{\L^c}) \right].
    \end{align*}
    
    For this, the strategy is similar to the classical proof of \cite[Theorem 15.37]{georgii_gibbs_2011} on the lattice. 
    Rougly speaking, the calculations above tell us that if $P \in \mathcal{P}_\Theta$ is a minimizer, then the relative entropy density $\delta F(P)$ vanishes. In particular, this implies that restrictions $P_\Lambda$ to sufficiently large but bounded volumes $\Lambda \subset \R^d$ are absolutely continuous with  density $g_{\Lambda}$ with respect to the finite-volume Gibbs measures $G_{\Lambda, z,\beta}$. 
    Moreover, the asymptotic control over the relative entropy in large volumes, which is implied by the vanishing relative entropy density, provides control over these densities for sufficiently large $\Lambda$. Making this precise is slightly technical but based on adapting standard arguments to our situation.

\section{Preliminaries}\label{sec_preliminaries}
    Let us start by collecting some properties of the conditional energy. Recall that we write $f\lesssim_a g$ in order to indicate that there exists a constant $c=c(a)$ that (only) depends to the parameter $a$, such that $f\le c g$.     
    \begin{lemma}[Properties of the energy]\label[lemma]{lem_energy}
        The energy function satisfies the following properties:
\begin{enumerate}
\item For all $\Lambda\in \mathfrak{B}_b(\R^d)$ we have that $H^{\rm hc}_{\Lambda}, H^{\rm ar}_{\Lambda}, H_{\Lambda}\ge 0$.
\item We have that
\begin{equation}\label{le_interaction_dominated_by_first_order}
H^{\rm ar}(\omega)
                \lesssim_d \sum_{\x \in \omega} (R_x+r)^d.
            \end{equation}
\item If \(\omega\in\Omega\) respects the hardcore condition, we have that
            \begin{equation}
                \sum_{\x \in \omega} R_x^d \lesssim_d H^{\mathrm{ar}}(\omega).
            \end{equation}
\item If \(\omega\in \Omega\) respects the hardcore condition and is finite, then
\begin{equation}\label{eq_x_wise_weighted_inclusion_exclusion}
H^{\mathrm{ar}}(\omega) 
                = \sum_{\x \in \omega} \Big(\phi_1(\x) + \sum_{k \geq 2} \tfrac{1}{k} \sum_{\{\y_1, \dots, \y_{k-1}\} \subseteq \omega\setminus\x} \phi_k(\x, \y_1, \dots, \y_{k-1}) \Big).
            \end{equation} 
            In particular,  for all \(\x \in \omega\) and \(\omega\) not necessarily finite,
\begin{equation}\label{eq_single_term_dominated_by_first_order}
\begin{split}
                \Big\vert \sum_{k \geq 2} &\tfrac{1}{k} \sum_{\{\y_1, \dots, \y_{k-1}\} \subseteq \omega\setminus \x} \phi_k(\x, \y_1, \dots, \y_{k-1}) \Big\vert\\
            &\leq \abs{B(x, R_x + r)\setminus B(x, R_x)} \lesssim_{d,r}(1+R_x^{d-1}).
            \end{split}
        \end{equation}
        
         \item For all $\zeta\in \Omega$ and \(\L\in \mathfrak{B}(\R^d)\), the map $\omega \mapsto H_{\L}(\omega_\L\zeta_{\L^c})$ is increasing.
        \end{enumerate}
    \end{lemma}
    \begin{proof}
        Note that Item~(1) is clear from the definitions.  
        Item~(2) follows from the estimates  $$
        H^{\mathrm{ar}}(\omega) = \Big\vert \bigcup_{\x \in \omega} B(x, R_x + r) \Big\vert \leq \sum_{\x \in \omega} \abs{B(x, R_x + r)} \lesssim_d \sum_{\x \in \omega} (R_x + r)^d.$$
For Item~(3), note that, by the hardcore condition, 
$$H^{\mathrm{ar}}(\omega) = \Big\vert \bigcup_{\x \in \omega} B(x, R_x + r)\Big\vert \geq \Big\vert \dot{\bigcup_{\x \in \omega}} B(x, R_x) \Big\vert = \sum_{\x \in \omega} \abs{B(x, R_x)} \gtrsim_d \sum_{\x \in \omega} R_x^d.$$
Now, for Item~(4), we have that for \(\x \in \omega\) it holds
        \begin{equation*}
        \begin{split}
            \frac{\sum_{\y \in \omega\setminus \x} \mathbbm{1}_{B(y, R_y + r)}}{\sum_{\y \in \omega} \mathbbm{1}_{B(y, R_y + r)}} &\mathbbm{1}_{B(x, R_x + r)} = \sum_{k \geq 2} \tfrac{(-1)^k}{k} \sum_{\{\y_1, \dots, \y_{k-1}\} \subseteq \omega\setminus \x} \mathbbm{1}_{B(x, R_x + r)} \prod_{i = 1}^{k-1} \mathbbm{1}_{B(y_i, R_{y_i} + r)}, 
            \end{split}
        \end{equation*}
        which can be seen by induction, 
        and therefore
        \begin{align*}
            &H^{\mathrm{ar}}(\omega)
            = \Big|\bigcup_{\x \in \omega} B(x, R_x + r)\Big|
            = \int \, \mathrm{d}z \, \mathbbm{1}\{z\in \bigcup_{\x \in \omega} B(x, R_x + r)\} \\
            &= \int \, \mathrm{d}z \, \sum_{\x \in \omega} \Big(\sum_{\y \in \omega} \mathbbm{1}\{z\in B(y, R_y + r)\}\Big)^{-1}  \mathbbm{1}\{z\in B(x, R_x + r)\}\\
            &= \int \, \mathrm{d}z \, \sum_{\x \in \omega} \Big(1 - \frac{\sum_{\y \in \omega\setminus \x} \mathbbm{1}\{z\in B(y, R_y + r)\}}{\sum_{\y \in \omega} \mathbbm{1}\{z\in B(y, R_y + r)\}} \Big) \mathbbm{1}\{z\in B(x, R_x + r)\} \\
            &= \sum_{\x \in \omega} \abs{B(x, R_x + r)} - \int \, \mathrm{d}z \,  \frac{\sum_{\y \in \omega\setminus \x} \mathbbm{1}\{z\in B(y, R_y + r)\}}{\sum_{\y \in \omega} \mathbbm{1}\{z\in B(y, R_y + r)\}} \mathbbm{1}\{z\in B(x, R_x + r)\} \\
            &= \sum_{\x \in \omega} \Big(\phi_1(x) + \sum_{k \geq 2} \tfrac{1}{k} \sum_{\{\y_1, \dots, \y_{k-1}\} \subseteq \omega\setminus \x} \phi_k(\x, \y_1, \dots, \y_{k-1}) \Big),
        \end{align*}
        as desired. Further, again by the hardcore condition,
        \begin{align*}
            &\Big\vert \sum_{k \geq 2} \tfrac{1}{k} \sum_{\{\y_1, \dots, \y_{k-1}\} \subseteq \omega\setminus\x} \phi_k(\x, \y_1, \dots, \y_{k-1})\Big\vert
            \\
            &= \int \, \mathrm{d}z \,  \frac{\sum_{\y \in \omega\setminus\x} \mathbbm{1}\{z\in B(y, R_y + r)\}}{\sum_{\y \in \omega} \mathbbm{1}\{z\in B(y, R_y + r)\}} \mathbbm{1}\{z\in B(x, R_x + r)\} \\
            &= \int \, \mathrm{d}z \,  \frac{\sum_{\y \in \omega\setminus\x} \mathbbm{1}\{z\in B(y, R_y + r)\}}{\sum_{\y \in \omega} \mathbbm{1}\{z\in B(y, R_y + r)\}} \mathbbm{1}\{z\in B(x, R_x + r)\setminus B(x, R_x)\}\\
            &\leq \abs{B(x, R_x + r)\setminus B(x, R_x)}
            \lesssim_{r,d} (1 + R_x^{d-1}).
        \end{align*}
        Note that, by the hardcore condition, only finitely many intersections of \(B(x, R_x +r)\) are non-empty such that everything above is well-defined even for a configuration \(\omega\) with infinitely many points.

        Item~(5) is again clear from the definition of \(H_{\L}(\omega_{\L}\zeta_{\L^c})\), as adding points to \(\omega\) clearly does not decrease \(H^{\mathrm{ar}}_{\L}(\omega_{\L}\zeta_{\L^c})\) or \(H^{\mathrm{hc}}_{\L}(\omega_{\L}\zeta_{\L^c})\).
    \end{proof}

Let us also collect some basic properties of the entropies. 
Here $\Fcal_\L$ denotes the sub-$\sigma$-algebra of events in $\Fcal$, which are measurable in $\L\in\mathfrak{B}(\R^d)$. 
    \begin{lemma}
        For probability measures \(\mu\) and \(\nu\) defined on the same probability space and a positive random variable \(f\) it holds that
        \begin{equation}
            \nu[f] \leq \log(\mu[\e^{f}]) + I(\nu \,\vert\, \mu).
        \end{equation}
        In particular, for \(P \in \mathcal{P}_{\Theta}\) and a positive \(\mathcal{F}_{\L}\)-measurable random variable \(f\) on \(\Omega\) with \( \L \in \mathfrak{B}_b(\mathbb{R}^d)\), we have
\begin{equation}\label{eq_entropic_inequality}
            P[f] 
            \leq \log\big(\pi^z[\e^{f}]\big) + \abs{\L} I^z(P).
        \end{equation}
    \end{lemma}
    \begin{proof}
        For the first inequality, see e.g.~\cite[Lemma 10.1]{Varadhan1984}.
        The second inequality is an application of the first to \(\nu = P_\L, \mu = \pi^z_\L\) combined with \Cref{eq_def_entropy_is_sup}.
    \end{proof}
\section{Energy density}\label{sec_energy_density}
Let us define, for any $x\in \R^d$, the spatial shift operator $\theta_x\colon \Omega\to\Omega$ with $\theta_x\omega=\theta_x((x_i,R_i)_{i\ge 0})=(x_i-x,R_i)_{i\ge 0}$. The following standard result establishes the Palm measure, that is the measure that describes a translation-invariant system from the perspective of a typical particle. Its proof can be found for example in~\cite[Remark 2.1]{GZ93}. 
    \begin{proposition}[Palm measure]\label[proposition]{prop_existence_palm_measure}
        For all \(P \in \mathcal{P}_{\Theta}\), there is a unique measure \(P^o\) on \(\left([0,\infty) \times \Omega, \mathfrak{B}([0,\infty)) \otimes \mathcal{F} \right)\), called the {\em Palm measure} of \(P\), such that
        \begin{equation}\label{eq_defining_property_palm_measure}
            \int_{\mathbb{R}^d} \mathrm{d}x \, \int_{\Omega \times [0,\infty)} P^o( \mathrm{d}\omega,\mathrm{d}R_o) \, f(x, R_o, \omega)
            = \int_{\Omega \times [0,\infty)} P(\mathrm{d}\omega) \, \sum_{\x \in \omega} f(x, R_x, \theta_x \omega)
        \end{equation} for all measurable \(f \colon \mathbb{R}^d \times [0,\infty) \times \Omega \to [0, \infty]\).
        In particular, \(P^o\) is a finite measure with \(P^o(\Omega \times [0,\infty)) = P[N_C] < \infty\), and
\begin{equation}\label{eq_radius_moment_palm_measure}
            P^o[R_o^d]
            = P\Bigl[\sum_{\x \in \omega_C} R_x^d \Bigr].
        \end{equation}
    \end{proposition}
    
Note that the Palm measure is supported on the set $\{(R,\omega)\in[0,\infty)\times\Omega \ \vert \ (o,R)\in\omega\}$.  
With Palm measures at hand, we can now state the following refined version of \Cref{prop_existence_of_energy_density}. 
\begin{proposition}[Existence and continuity properties of the energy density]\label[proposition]{prop_existence_of_energy_density_2}
For all \(P \in \mathcal{P}_{\Theta}\), the energy density 
        \begin{align*}
            H(P) := \lim_{n \uparrow \infty} \frac{P[H_{\L_n, \emptyset}]}{\abs{\L_n}}
        \end{align*} exists and is given by
        \begin{align*}
            H(P) = \begin{cases}
                P^o[f_H], & \text{if } P\Big[\sum_{\{\x, \y \} \subseteq \omega} \phi_{\mathrm{hc}}(\x, \y)\Big] < \infty \text{ and } P\Big[\sum_{\x \in\omega_{C}} R_x^d\Big] < \infty, \\
                +\infty, & \text{ otherwise},
            \end{cases} 
        \end{align*} where \(P^o\) is the Palm measure of \(P\), and
        \begin{align*}
            f_H(R_o, \omega) 
            := 
            \phi_1(o, R_o) + \sum_{k \geq 2} \tfrac{1}{k} \sum_{\{\x_1, \dots, \x_{k-1}\} \subseteq \omega\setminus \oo} \phi_k((o, R_o), \x_1, \dots, \x_{k-1}).
        \end{align*}
        Additionally, \(P \mapsto H(P)\) is continuous on \(\{P\in P_\Theta\,\vert\, H(P) < \infty\} \) and in particular lower semicontinuous on \(P_{\Theta}\) w.r.t.~the \(\tau_\mathcal{L}\)-topology.
    \end{proposition}
    The proof of \Cref{prop_existence_of_energy_density_2} is based on a few technical lemmas, which we provide now. The proof of \Cref{prop_existence_of_energy_density_2} is given at the end of this section. 
\begin{lemma}[Finite-volume Palm representation]\label[lemma]{le:representation_of_energy_in_finite_volume}
        Let \(P \in \mathcal{P}_{\Theta}\). If
        \begin{align*}
            P\Big[\sum_{\{\x, \y \} \subseteq \omega} \phi_{\mathrm{hc}}(\x, \y)\Big] < \infty \quad\text{and}\quad P\Big[\sum_{\x \in\omega_{C}} R_x^d\Big] < \infty,
        \end{align*}
        then
        \begin{align*}
            \frac{P[H_{\L_n, \emptyset}]}{\abs{\L_n}}
            = P^o[f_{H, n}],
        \end{align*} 
        where
        \begin{align*}
            f_{H, n}&(R_o, \omega):= \phi_1(\oo)+ \sum_{k \geq 2} \tfrac{1}{k} \sum_{\{\x_1, \dots, \x_{k-1}\} \subseteq \omega\setminus \oo} \phi_k(\oo, \x_1, \dots, \x_{k-1}) \tfrac{\big|\bigcap_{i = 1}^{k-1} (\L_n - x_i) \cap \L_n\big|}{\abs{\L_n}}.
        \end{align*}
    \end{lemma}
    \begin{proof}
        By \Cref{eq_x_wise_weighted_inclusion_exclusion}, we have that
        \begin{align*}
            H^{\mathrm{ar}}_{\L_n, \emptyset}(\omega)
            &= \Big|\bigcup_{\x \in \omega_{\L_n}} B(x, R_x + r)\Big|\\
            &= \sum_{\x \in \omega_{\L_n}} \Big(\phi_1(\x) + \sum_{k \geq 2} \tfrac{1}{k} \sum_{\{\y_1, \dots, \y_{k-1}\} \subseteq \omega_{\L_n}\setminus\x} \phi_k(\x, \y_1, \dots, \y_{k-1}) \Big).
        \end{align*}
        Hence, by translation invariance, and since under \(P\) we can ignore the hardcore term,
        \begin{align*}
            H_{\L_n, \emptyset}(\omega)
            = H^{\mathrm{ar}, \emptyset}_{\L_n}(\omega)
            = \sum_{\x \in \omega} g_n(\x, \theta_x \omega),
        \end{align*} where
        \begin{align*}
            g_n(\x, \omega)
            := \phi_1&((o, R_x)) \mathbbm{1}\{x\in \L_n\}\\
            &\hspace{-1cm}+ \sum_{k \geq 2} \tfrac{1}{k} \sum_{\{\y_1, \dots, \y_{k-1}\} \subseteq \omega\setminus (o,R_x)}\hspace{-.8cm} \phi_k((o, R_x), \y_1, \dots, \y_{k-1}) \mathbbm{1}\Big\{x\in \bigcap_{i=1}^{k-1} (\L_n - y_i) \cap \L_n\Big\}.
        \end{align*}
        The defining property of the Palm measure \(P^o\), \Cref{eq_defining_property_palm_measure}, gives
        \begin{align*}
            P[H_{\L_n, \emptyset}]
            &= P\Big[\sum_{\x \in \omega} g_n(\x, \theta_x \omega)\Big]
            = \int \mathrm{d}x \int P^o( \mathrm{d}\omega,\mathrm{d}R_o) \, g_n((x, R_o), \omega) \\
            &= \int P^o( \mathrm{d}\omega,\mathrm{d}R_o) \int \mathrm{d}x \, g_n((x, R_o), \omega) \\
            &= P^o[\abs{\L_n} f_{H, n}],
        \end{align*}
        as desired. 
    \end{proof}
\begin{lemma}[Infinite-volume Palm representation]\label[lemma]{le:representation_of_energy_in_infinite_volume}
        Let \(P \in \mathcal{P}_{\Theta}\) be such that both \(P\big[\sum_{\{\x, \y \} \subseteq \omega} \phi_{\mathrm{hc}}(\x, \y)\big] < \infty\) and \(P\big[\sum_{\x \in\omega_{C}} R_x^d\big] < \infty\).
        Then \(P^o(f_H)\) is well-defined for
        \begin{align*}
            f_H(R_o, \omega) 
            := \phi_1(\oo) + \sum_{k \geq 2} \tfrac{1}{k} \sum_{\{\x_1, \dots, \x_{k-1}\} \subseteq \omega\setminus \oo} \phi_k(\oo, \x_1, \dots, \x_{k-1})
        \end{align*}
        and furthermore, $
           P^o(f_{H}) = \lim_{n \uparrow \infty} P^o[f_{H, n}]$.
    \end{lemma}
    \begin{proof}
        The random variable \(f_H\) is integrable w.r.t.~\(P^o\), since \(\phi_1\) is integrable by the moment assumption on \(P\) (see \Cref{eq_radius_moment_palm_measure}) and, by \Cref{eq_single_term_dominated_by_first_order},
        \begin{align*}
            \left\vert f_H(R_o, \omega) 
            - \phi_1(\oo) \right\vert
            \lesssim \phi_1(\oo).
        \end{align*}
        The dominated convergence theorem now implies \(P^o(f_H) = \lim_{n \uparrow \infty} P^o(f_{H, n})\).
    \end{proof}
    \begin{proposition}[Continuity of the energy density]\label[proposition]{thm_lsc_of_energy_density}
        The map \(P \mapsto H(P)\) described in \Cref{prop_existence_of_energy_density_2} is continuous on \(\{P\in \mathcal P_\Theta\colon H(P) < \infty\}\) and in particular lower-semicontinuous w.r.t.~the \(\tau_\mathcal{L}\)-topology on $\mathcal{P}_\Theta$.
    \end{proposition}
    \begin{proof}
        Let \((P_\alpha)_\alpha\) be any net in \(\{P\in \mathcal P\colon H(P) < \infty\}\) which converges to \(P \in \mathcal{P}_{\Theta}\) in the \(\tau_\mathcal{L}\)-topology.
        We first establish that \(P\) has to respect the hardcore condition. Clearly, every \(P_\alpha\) has to respect the hardcore condition, since \(H(P_\alpha) < \infty\).
        Moreover, since
        \begin{align*}
            P\Big(\sum_{\{\x, \y\} \subseteq \omega} \phi_{\mathrm{hc}}(\x, \y) = \infty \Big)
            = \sup_{n \in \mathbb{N}} P\Big(\sum_{\{\x, \y\} \subseteq \omega_{\Lambda_n}} \, \phi_{\mathrm{hc}}(\x, \y) = \infty \Big)
        \end{align*} 
        and by the definition of the \(\tau_\mathcal{L}\)-topology
        \begin{align*}
            0
            = P_\alpha\Big(\sum_{\{\x, \y\} \subseteq \omega_{\Lambda_n}} \, \phi_{\mathrm{hc}}(\x, \y) = \infty \Big)
            \to P\Big(\sum_{\{\x, \y\} \subseteq \omega_{\Lambda_n}} \, \phi_{\mathrm{hc}}(\x, \y) = \infty \Big),
        \end{align*} 
        we get
        \begin{align*}
            P\Big(\sum_{\{\x, \y\} \subseteq \omega} \phi_{\mathrm{hc}}(\x, \y) = \infty\Big) = 0.
        \end{align*}
        By the definition of the \(\tau_\mathcal{L}\)-topology, we also have that
        \begin{align*}
             P_\alpha\Big[\sum_{\x \in \omega_C} R_x^d \Big]
            \to P\Big[\sum_{\x \in \omega_C} R_x^d \Big],
        \end{align*}
        so that, either 
        \begin{align*}
            H(P) \gtrsim 
            P\Big[\sum_{\x \in \omega_C} R_x^d \Big]
            = \infty
        \end{align*} 
        and also \(\lim_{\alpha } H(P_\alpha) \gtrsim \lim_{\alpha} P_\alpha[\sum_{\x \in \omega_c} R_x^d] = \infty\) by \Cref{le_interaction_dominated_by_first_order}, or
        \begin{align*}
            P\Big[\sum_{\x \in \omega_C} R_x^d\Big] =
            \lim_{\alpha} P_\alpha\Big[\sum_{\x \in \omega_C} R_x^d\Big]
            \leq c
        \end{align*} for some \(c > 0\), and we can assume without loss of generality that
        \begin{align*}
            \sup_{\alpha} P_\alpha^o[R_o^d]
            = \sup_{\alpha} P_\alpha\Big[\sum_{\x \in \omega_C} R_x^d\Big]
            \lesssim \sup_{\alpha} H(P_\alpha)
            \leq c,
        \end{align*} where we used \Cref{eq_radius_moment_palm_measure} and \Cref{le_interaction_dominated_by_first_order}.
        Given that this holds, we have
        \begin{align*}
            H(P)
            = P^o\left[\phi_1(\oo)
            + X\right],
        \end{align*} where
        \begin{align*}
            X := 
            X(R_o, \omega) 
            := \sum_{k \geq 2} \tfrac{1}{k} \sum_{\{\x_1, \dots, \x_{k-1}\} \subseteq \omega\setminus\oo} \phi_k(\oo, \x_1, \dots, \x_{k-1}).
        \end{align*}
         We now show that $            P_\alpha^o\left[X\right]
            \to
            P^o\left[X\right]$.         Fix an \(\epsilon > 0\) and consider boxes \(\Delta_1, \dots, \Delta_M \subseteq \mathbb{R}^d\) to be specified later.
        Then, for \(Q \in \{P, P_\alpha\}\), we have
        \begin{align*}
            &\Big\vert Q^o\left[X \right]
            -
            Q^o\big[X \mathbbm{1}_{N_{\Delta_1}, \dots, N_{\Delta_M} > 0 \text{ and } R_o \leq \epsilon^{-1}} \big]  \Big\vert
            \leq \sum_{m = 1}^{M} Q^o\left[\abs{X} \mathbbm{1}_{N_{\Delta_m} = 0} \right]
            + Q^o\left[\abs{X} \mathbbm{1}_{R_o > \epsilon^{-1}} \right].
        \end{align*}
        Recall that \(\abs{X} \lesssim 1 + R_o^{d-1} \lesssim R_o^{d-1}\) for \(R_o \geq 1\) by \Cref{eq_single_term_dominated_by_first_order}. Hence,
        \begin{align*}
            Q^o\big[\abs{X}  \mathbbm{1}_{R_o > \epsilon^{-1}} \big]
            &\leq \epsilon Q^o\big[R_o \abs{X} \mathbbm{1}_{R_o > \epsilon^{-1}}  \big]
            \lesssim \epsilon Q^o\big[R_o R_o^{d-1}  \big]
            \lesssim \epsilon.
        \end{align*}
        We also have
        \begin{align*}
            &Q^o\left[\abs{X}  \mathbbm{1}_{N_{\Delta} = 0} \right] 
            = Q\Big[\sum_{\x \in \omega_C} \abs{X(R_x, \theta_x \omega)} \mathbbm{1}_{N_{\Delta} = 0}(\theta_x \omega) \Big]\\
            &\leq  Q\Big[\Big(\sum_{\x \in \omega_C} \abs{X(R_x, \theta_x \omega)}\Big) \mathbbm{1}_{N_{\Delta^{-}} = 0}(\omega) \Big]\\
            &= Q\Big[\Big(\sum_{\x \in \omega_C} \Big\vert \sum_{k \geq 2} \tfrac{1}{k} \sum_{\{\x_1, \dots, \x_{k-1}\} \subseteq \omega\setminus \x} \phi_k(\x, \x_1, \dots, \x_{k-1}) \Big\vert\Big) \mathbbm{1}_{N_{\Delta^{-}} = 0}(\omega) \Big] \\
            &\leq Q\Big[\Big(\sum_{\x \in \omega_C} \vert B(x, R_x + r) \setminus B(x, R_x) \vert\Big) \mathbbm{1}_{N_{\Delta^{-}} = 0, N_C \neq 0}(\omega) \Big]\\
            &\leq Q\Big[\Big(\sum_{\x \in \omega_C} \vert B(x, R_x + r) \setminus B(x, R_x) \vert\Big)^{d/(d-1)} \Big]^{(d-1)/d} Q(N_{\Delta^{-}} = 0, N_C \neq 0)^{1/d},
        \end{align*} 
        where \(\Delta^{-}\) is the maximal box inside \(\Delta\) such that \(\theta_{x} \Delta^{-} \subseteq \Delta\) for all \(x \in C = \Lambda_1\), and we used \Cref{eq_single_term_dominated_by_first_order}.
        However, by the hardcore condition,
        \begin{align*}
            &\Big(\sum_{\x \in \omega_C} \vert B(x, R_x + r) \setminus B(x, R_x) \vert\Big)^{d/(d-1)}\\
            &= 
            \Big(\sum_{\x \in \omega_C} \vert B(x, R_x + r) \setminus B(x, R_x) \vert\Big)^{d/(d-1)} \mathbbm{1}\{N_C > 1\} \\
            &\qquad+ \Big(\sum_{\x \in \omega_C} \vert B(x, R_x + r) \setminus B(x, R_x) \vert\Big)^{d/(d-1)} \mathbbm{1}\{N_C = 1\}\\
            &\lesssim \vert C \vert^{d/(d-1)} + 1 + \Big(\sum_{\x \in \omega_C} R_x^d \Big)
            \lesssim 1 + \sum_{\x \in \omega_C} R_x^d,
        \end{align*} 
        and hence
        \begin{align*}
            Q^o\Big[\abs{X}  \mathbbm{1}_{N_{\Delta} = 0} \Big]
            &\lesssim Q\Big[1 + \sum_{\x \in \omega_C} R_x^d \Big]^{(d-1)/d} Q(N_{\Delta^{-}} = 0, N_C \neq 0)^{1/d}\\
            &= \Big(1 + Q^o[R_o^d] \Big)^{(d-1)/d} Q(N_{\Delta^{-}} = 0, N_C \neq 0)^{1/d}\\
            &\lesssim Q(N_{\Delta^{-}} = 0, N_C \neq 0)^{1/d}.
        \end{align*}
        For a fixed box \(\Delta\), we then have
        \begin{align}
             P_\alpha(N_{\Delta^{-}} = 0, N_C \neq 0)
             \to 
             P(N_{\Delta^{-}} = 0, N_C \neq 0)
        \end{align} by the definition of the \(\tau_\mathcal{L}\)-topology.
        On the other hand, for \(\Delta \uparrow \mathbb{R}^d\) we have
        \begin{align*}
            P(N_{\Delta^{-}} = 0, N_C \neq 0)
            &\leq P(N_{\Delta^{-}} = 0, \omega \neq \emptyset)
            = P(N_{\Delta_o^{-}} = 0, \omega \neq \emptyset)\\
            &\downarrow
            P(N_{\mathbb{R}^d} = 0, \omega \neq \emptyset)
            = 0,
        \end{align*} 
        by translation invariance, where \(\Delta_o^-\) is the translate of the box \(\Delta^{-}\) which is centered at the origin.
    
        In summary, if we choose \(\Delta\) large enough, we have
        \begin{align*}
           P_\alpha(N_{\Delta^{-}} = 0, N_C \neq 0)
            \to 
             P(N_{\Delta^{-}} = 0, N_C \neq 0)
            < \epsilon
        \end{align*} and can therefore assume without loss of generality that
        \begin{align*}
            \left\vert Q^o\left[X \right]
            -
            Q^o\left[X \mathbbm{1}_{N_{\Delta_1}, \dots, N_{\Delta_M} > 0 \text{ and } R_o \leq \epsilon^{-1}} \right]  \right\vert
            \lesssim M \epsilon^{\frac{1}{d}}
            + \epsilon
        \end{align*} for \(Q \in \{P, P_\alpha\}\).
        Therefore, we only need to provide some fixed number \(M \in \mathbb{N}\) such that there are arbitrarily large sets \(\Delta_1, \dots, \Delta_M\) with
        \begin{align*}
            P_\alpha^o\big[X \mathbbm{1}_{N_{\Delta_1}, \dots, N_{\Delta_M} > 0 \text{ and } R_o \leq \epsilon^{-1}}\big]  
            \to
            P^o\big[X \mathbbm{1}_{N_{\Delta_1}, \dots, N_{\Delta_M} > 0 \text{ and } R_o \leq \epsilon^{-1}}\big]. 
        \end{align*}
        \begin{figure}[t]
             \centering 
            \includegraphics[width=.6\textwidth]{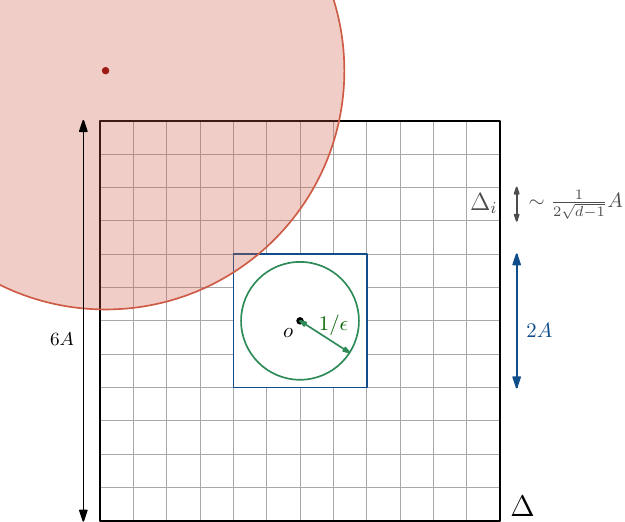}
            \caption{Such a geometric arrangement forces balls with centers outside \(\Delta\), which intersect \([-A, A]^d \supset B(0,\epsilon^{-1})\), to contain at least one of the \(\Delta_j\).}
            \label{fig:geometric_blocking}
        \end{figure}
        We will do this by choosing \(\Delta_1, \dots, \Delta_M\) in such a manner that there is a bounded \(\Delta\) with
        \begin{equation}\label{eq:interactions_can_be_localized}
        \begin{split}
            &X \mathbbm{1}_{\{N_{\Delta_1}, \dots, N_{\Delta_M} > 0 \text{ and } R_o \leq \epsilon^{-1}\}}\\
            &= \Big(\sum_{k \geq 2} \tfrac{1}{k} \sum_{\{\x_1, \dots, \x_{k-1}\} \subseteq \omega_{\Delta}\setminus\oo} \phi_k(\oo, \x_1, \dots, \x_{k-1})\Big) \mathbbm{1}_{\{N_{\Delta_1}, \dots, N_{\Delta_M} > 0 \text{ and } R_o \leq \epsilon^{-1}\}},
        \end{split}
        \end{equation}
        i.e., such that we actually only have to look at interactions in the fixed, bounded volume \(\Delta\). This can be achieved by choosing, for some scale parameter \(A \geq \epsilon^{-1}\), the set \(\Delta = [-3A, 3A]^d\) and decomposing the ``boundary'' \(\Delta \setminus [-A, A]^d\) into \(M \sim 2^d (3^d-1) (d-1)^{d/2}\) rectangles with edge length \(\sim \frac{A}{2\sqrt{d-1}}\). See \Cref{fig:geometric_blocking} for a sketch. This geometric arrangement now has the property that every ball with center outside \(\Delta\) which intersects \([-A, A]^d \supset B(0,\epsilon^{-1})\) has to fully contain \(\Delta_j\) for at least one \(j \in \{1, \dots, M\}\). The latter property together with the hardcore condition in our model guarantees \Cref{eq:interactions_can_be_localized}.
    \end{proof}
We are now in the position to prove the main result of this section. 
\begin{proof}[Proof of \Cref{prop_existence_of_energy_density_2}]
        If \(P\left[\sum_{\{\x, \y \} \subseteq \omega} \phi_{\mathrm{hc}}(\x, \y)\right] = \infty\) or \(P\left[\sum_{\x \in\omega_{C}} R_x^d\right] = \infty\), then \(P[H_{\L_n}] = \infty\) for all \(n \in \mathbb{N}\) sufficiently large, by~\Cref{le_interaction_dominated_by_first_order}.
        Therefore, we can assume \(P\big[\sum_{\{\x, \y \} \subseteq \omega} \phi_{\mathrm{hc}}(\x, \y)\big] < \infty\) and \(P\big[\sum_{\x \in\omega_{C}} R_x^d\big] < \infty\) from now on.
        Then, we have \(P[H_{\L_n, \emptyset}]/\abs{\L_n} = P^o[f_{H, n}]\) by \Cref{le:representation_of_energy_in_finite_volume} and therefore we have that
        \begin{align*}
            H(P)
            = \lim_{n \uparrow \infty}\frac{P[H_{\L_n, \emptyset}]}{\abs{\L_n}} 
            = \lim_{n \uparrow \infty} P^o[f_{H, n}]
            =  P^o[f_{H}]
        \end{align*} by \Cref{le:representation_of_energy_in_infinite_volume}.
        The continuity properties of \(H\) are shown in \Cref{thm_lsc_of_energy_density}.
    \end{proof}
    
    Before the start of the next section, let us note the following.  
    \begin{remark}\label{remark_discontinuity}
        If infinite energies are allowed, the map \(P \mapsto H(P)\) is not continuous. 
        Here is an example of translation-invariant probability measures \((P_n)_{n \in \mathbb{N}}\) such that \(H(P_n) = \infty\) for all \(n \in \mathbb{N}\) and \(P_n \to P\) in the \(\tau_\mathcal{L}\)-topology, but \(H(P) < \infty\).
        Define the periodic hard-core-respecting configuration
        \begin{align*}
            \omega_{\text{good}}
            = \sum_{i \in \mathbb{Z}^d} \delta_{(i, 1/3)},
        \end{align*}
        and the periodic non-hard-core-respecting (i.e., overlapping) configuration
        \begin{align*}
            \omega_{\text{bad}}
            = \sum_{i \in \mathbb{Z}^d} \delta_{(i, 2)}.
        \end{align*} 
        Finally, define the translation-invariant states
        \begin{align*}
            \mu_{\text{good}} 
            = 2^{-d} \int_{[-1, 1]^d} \, \delta_{\theta_y  \omega_{\text{good}}} \, \mathrm{dy}\qquad\text{ and }\qquad
            \mu_{\text{bad}} 
            = {2}^{-d} \int_{[-1, 1]^d} \, \delta_{\theta_y \omega_{\text{bad}}} \, \mathrm{dy},
        \end{align*}
        and set $P_n = (1 - 1/n) \mu_{\text{good}} + (1/n) \mu_{\text{bad}}$. 
Then, we have that $H(P_n) \geq  (1/n)H(\mu_{\text{bad}}) = \infty$, while on the other hand $H(\mu_{\text{good}}) < \infty$ and $P_n \xrightarrow[n\to \infty]{} \mu_{\text{good}} =: P$ in the \(\tau_{\mathcal{L}}\)-topology.
    \end{remark}

\section{The support of good measures: Temperedness}\label{sec_temperedness}

Let us start by defining the following subsets of configurations related to the temperedness condition. For $K,N,M \in \N$ with $N \leq M$ we will write
        \begin{align*}
            \Omega^*_{K}
            &:= \left\{\omega \in \Omega \ \middle\vert \ \forall n \geq K \ \forall \x \in \omega_{\L_n} \colon  R_x \leq g(n)  \right\},
             \\
            \Omega^*_{N, M} &:= \left\{\omega \in \Omega \ \middle\vert \ \forall n \in \mathbb{N} \cap [N, M] \ \forall \x \in \omega_{\L_n} \colon  R_x \leq g(n)  \right\}.
        \end{align*} 
        For \(\L \in \mathfrak{B}(\mathbb{R}^d)\), we let
\begin{align}\label{sets_of_temperedness_conditions}
            (\Omega_{\L})^{*}_{K} 
            = \left\{\omega \in \Omega_\L \ \middle\vert \ \forall n \geq K \ \forall \x \in \omega_{\L_n \cap \L} \colon R_x \leq g(n)  \right\}
        \end{align}
        and we will later also need
        \begin{align*}
            \Omega^{*, \L}_{K}  
            = \left\{\omega \in \Omega \ \middle\vert \ \forall n \geq K \ \forall \x \in \omega_{\L_n \cap \L} \colon R_x \leq g(n)  \right\}.
        \end{align*}

    With this notation at hand, we first derive some quantitative estimates for the underlying Poisson point process.  

    \begin{lemma}[Temperedness estimates in the Poisson case]\label[lemma]{le_temperedness_estimate_poisson}
       We have that 
        \begin{align*}
            \pi^z((\Omega^*_{N, M})^c)
            \lesssim_{\mathfrak{R}, \gamma, z} \e^{-\abs{\L_N}^{1+\gamma}}
        \end{align*} 
        and in particular, for all \(K \in \mathbb{N}\),
        \begin{align*}
\pi^z((\Omega^*_{K})^c)
            = \sup_{M \geq K} \pi^z((\Omega^*_{K, M})^c)
            \lesssim_{\mathfrak{R}, \gamma,z} \e^{-\abs{\L_N}^{1+\gamma}}.
        \end{align*} 
    \end{lemma}
    \begin{proof}
We have that
        \begin{align*}
            \pi^z(\Omega^*_{N, M})
            &= \pi^z(\forall \x \in \omega_{\L_N} \colon R_x \leq g(N)) \prod_{n = N+1}^{M} \pi^z(\forall \x \in \omega_{\L_n \setminus \L_{n-1}} \colon R_x \leq g(n)) \\
            &= \exp\Big(-z\abs{\L_N}\mathfrak{R}(R > g(N)) - z\sum_{n = N+1}^{M} \abs{\L_{n}\setminus \L_{n-1}}\mathfrak{R}(R > g(n)) \Big),
        \end{align*} which implies that
        \begin{align*}
            \pi^z((\Omega^*_{N, M})^c)
            &= 1 - \exp\Big(-z\abs{\L_N}\mathfrak{R}(R > g(N)) - z\sum_{n = N+1}^{M} \abs{\L_{n}\setminus \L_{n-1}}\mathfrak{R}\big(R > g(n)\big) \Big) \\
            &\leq 1 - \exp\Big(-z\abs{\L_N}\mathfrak{R}(R > g(N)) - z\sum_{n = N+1}^{\infty} \abs{\L_{n}\setminus \L_{n-1}}\mathfrak{R}\big(R > g(n)\big) \Big) \\
            &\leq 1 - \exp\Big(- c'_{\mathfrak{R}, \gamma,z}\e^{-\abs{\L_N}^{1+\gamma}} - c''_{\mathfrak{R}, \gamma,z}\e^{-\abs{\L_N}^{1+\gamma}}\Big)
            \lesssim_{\mathfrak{R}, \gamma,z} \exp(-\abs{\L_N}^{1+\gamma}),
        \end{align*} 
        as desired. 
    \end{proof}

    We directly apply this technical helper to show that measures $P$ with finite specific entropy are necessarily tempered. 
    \begin{proof}[Proof of \Cref{lemma_support_of_measure_with_finite_free_energy}]
    Let us note that the obvious union bound does not work for \(d=1\) because then \(\abs{\L_n}^{-\gamma}\) is not summable for \(\gamma \in (0,1)\), hence the following slightly more complicated analysis. 
    From the general inequality \Cref{eq_entropic_inequality} and \Cref{le_temperedness_estimate_poisson} we have for all \(N, M \in \mathbb{N}\) with \(N \leq M\),
        \begin{align*}
            P\big((\Omega^{*}_{N, M})^c\big)
            &= \frac{P[\abs{\L_M}^{1+\gamma/2} \mathbbm{1}\{(\Omega^{*}_{N, M})^c\}]}{\abs{\L_M}^{1+\gamma/2}}\\
            &\leq \frac{\log\big( \e^{\abs{\L_M}^{1+\gamma/2}}\pi^z\big((\Omega^{*}_{N, M})^c\big) + 1 \big)}{\abs{\L_M}^{1+\gamma/2}} + \frac{I^z(P)}{\abs{\L_M}^{\gamma/2}}\\
            &\leq \frac{\log\big(c_{\mathfrak{R}, \gamma,z} \e^{\abs{\L_M}^{1+\gamma/2}} \e^{-\abs{\L_N}^{1+\gamma}} + 1 \big)}{\abs{\L_M}^{1+\gamma/2}} + \frac{I^z(P)}{\abs{\L_M}^{\gamma/2}}.
        \end{align*}
        For the special choice \(N = N_k\) and \(M = N_{k+1}\) where \(N_k = k^{4/\gamma}\) we thus have that
        \begin{align*}
            P((\Omega^{*}_{N_k, N_{k+1}})^c)
            \lesssim \frac{\log\Big(c_{\mathfrak{R}, \gamma,z} \e^{-c_d k^{4(1+1/\gamma)}(1 - k^{-2}((k+1)/k)^{4(1/2 + 1/\gamma)})} + 1 \big)}{k^{(4+1/2)d}} + \frac{I^z(P)}{k^{2d}}
            \lesssim k^{-2d}.
        \end{align*}
        Since \(\Omega^* = \bigcup_{K \in \mathbb{N}} \bigcap_{k \geq K} \Omega^*_{N_k, N_{k+1}}\) we finally have
        \begin{align*}
            P\big((\Omega^{*})^c\big)
            \leq \inf_{K \in \mathbb{N}} \sum_{k \geq K}  P\big((\Omega^{*}_{N_k, N_{k+1}})^c\big)
            \lesssim \inf_{K \in \mathbb{N}} \sum_{k \geq K} k^{-2d} 
            = 0,
        \end{align*}
        as desired. 
    \end{proof}
    
    In order to show that infinite-volume translation-invariant Gibbs measures are tempered, we first provide a finite-volume comparison result. 
    
    \begin{lemma}[Temperedness estimates for finite-volume Gibbs measures]\label[lemma]{le_temperedness_estimate_gibbs}
        Let \(\L \in \mathfrak{B}_b(\mathbb{R}^d)\) and \(K \in \mathbb{N}\), then, for every boundary configuration \(\zeta \in \Omega\),
        \begin{align*}
            G_{\L, z, \beta, \zeta}\big((\Omega_\L)^{*}_{K}\big)
            \geq \pi^z(\Omega^{*}_{K}). 
        \end{align*} 
    \end{lemma}
    \begin{proof}
        Since \(\omega_\L \mapsto \exp(-\beta H_{\L, \zeta}(\omega_\L))\) and \(\omega_\L \mapsto \mathbbm{1}\{\omega_\L\in (\Omega_\L)^{*}_{K}\}\) are both decreasing, the result follows by applying the Harris-/FKG-inequality for the Poisson point process $\pi^z_{\L}$, i.e.,
        \begin{align*}
            G_{\L, z, \beta, \zeta}\big((\Omega_\L)^{*}_{K}\big)
            &= \pi^z_{\L}\left[\frac{\e^{-\beta H_{\L, \zeta}}}{Z_{\L, z, \beta, \zeta}} \mathbbm{1}_{(\Omega_\L)^{*}_{K}}\right]
            \geq \pi^z_{\L}\left[\frac{\e^{-\beta H_{\L, \zeta}}}{Z_{\L, z, \beta, \zeta}} \right] \pi^z_{\L}((\Omega_\L)^{*}_{K}) \\
            &=\pi^z_{\L}((\Omega_\L)^{*}_{K})
            \geq \pi^z(\Omega^{*}_{K}).
        \end{align*}
    \end{proof}

    This finite-volume comparison now allows us to conclude that the translation-invariant infinite-volume Gibbs measures are also tempered. 
    
    \begin{proof}[Proof of  \Cref{le_infinite_volume_gibbs_measures_are_tempered}]
        By the definition of infinite-volume Gibbs measures and \Cref{le_temperedness_estimate_gibbs}, we have that
        \begin{align*}
            P(\Omega_K^{*, \L})
           & = \int P(\mathrm{d}\zeta) \int G_{\L, z, \beta, \zeta}(\mathrm{d}\omega_\L) \mathbbm{1}\{\omega_\L \zeta_{L^c}\in \Omega_K^{*, \L}\}
            = \int P(\mathrm{d}\zeta) G_{\L, z, \beta, \zeta}\big((\Omega_\L)_K^{*}\big)\\
            &\geq \int P(\mathrm{d}\zeta) \pi^z(\Omega^{*}_{K})
            = \pi^z(\Omega^{*}_{K}).
        \end{align*}
        But then, using \Cref{le_temperedness_estimate_poisson},
        \begin{align*}
            P\left(\Omega^*\right)
            = \sup_{K \in \mathbb{N}} P\left(\Omega_K^*\right)
            = \sup_{K \in \mathbb{N}} \inf_{\Lambda \in \mathfrak{B}_b(\mathbb{R}^d)} P(\Omega_K^{*, \L})
            \geq \sup_{K \in \mathbb{N}} \pi^z(\Omega^{*}_{K})
            = 1.
        \end{align*}
    \end{proof}
    We conclude this section with another regularity result for the concatenation of tempered measures with the finite-volume Gibbs specification. This will be important for the proof of \Cref{thm_variational_principle} in \Cref{sec_minimizers_are_gibbs}. 
    \begin{lemma}\label[lemma]{le_temperedness_estimate_for_nu}
        Let \(P \in \mathcal{P}_{\Theta}\) be tempered.
        Then, for all \(\L \in \mathfrak{B}_b(\mathbb{R}^d)\), the measure \(\nu_\L = PG_{\L}\) defined, for any $A \in \mathcal{F}$, by
        \begin{align*}
            \nu_\L(A) = \int P(\mathrm{d}\zeta) \int  G_{\L, z, \beta, \zeta}(\mathrm{d}\omega_\L) \, \mathbbm{1}\{\omega_{\L} \zeta_{\L^c}\in A\}
        \end{align*} 
        is also tempered.
        Moreover,
        \begin{align*}
            \sup_{\L \in \mathfrak{B}_b(\mathbb{R}^d)} \nu_\L({\Omega^*_{K}}^c)
            \xrightarrow[K \uparrow \infty]{} 0.
        \end{align*}
    \end{lemma}
    \begin{proof}
        For any $\L \in \mathfrak{B}_b(\mathbb{R}^d)$, using \Cref{le_temperedness_estimate_gibbs}, we compute
        \begin{align*}
            \nu_\L(\Omega^*_{K})
            &= PG_\L(\Omega^*_{K})
            = \int \, P(\mathrm{d}\zeta) \, \int \, G_{\L, z, \beta, \zeta}(\mathrm{d}\omega) \mathbbm{1}\{\omega_{\L}\zeta_{\L^c}\in \Omega^*_{K}\}\\
            &\geq \int \, P(\mathrm{d}\zeta) \, \int \, G_{\L, z, \beta, \zeta}(\mathrm{d}\omega) \mathbbm{1}\{\zeta\in \Omega^{*, \L}_{K}\} \mathbbm{1}\{\omega_{\L}\in (\Omega_\L)^{*}_{K}\}\\
            &= \int \, P(\mathrm{d}\zeta) \, \mathbbm{1}\{\zeta\in \Omega^*_{K}\}  \, G_{\L, z, \beta, \zeta}\big((\Omega_\L)^{*}_{K}\big)\\
            &
            \geq \int \, P(\mathrm{d}\zeta) \, \mathbbm{1}\{\zeta\in \Omega^*_{K}\}  \, \pi^z(\Omega^{*}_{K})\\
            &= P(\Omega^*_{K}) \pi^z(\Omega^{*}_{K}),
        \end{align*} 
        where the right-hand side is independent of \(\L\) and tends to one as $K\uparrow\infty$ by the temperedness assumption on \(P\) and \Cref{le_temperedness_estimate_poisson}.
\end{proof}

\section{Existence of pressure via large deviations}\label{sec_pressure}
In this section, we show the existence of the density limit of the partition function. We do this step-by-step by exchanging boundary conditions. We start by considering the periodic case and then extend to free and tempered boundary conditions via comparison estimates. 

\subsection{Periodic boundary conditions}
    The fundamental idea is to use an underlying large-deviation result for the stationary empirical field
    \begin{align*}
        R_{n, \omega} := \abs{\L_n}^{-1} \int_{\Lambda_n} \, \delta_{\theta_x \omega^{(n)}} \, \mathrm{d}x.
    \end{align*}
    
    Recall that  $H_{\L_n, \mathrm{per}}(\omega_{\L_n})$ denotes the finite-volume energy when $\L_n$ is identified with the $d$-dimensional torus. 
    The precise relation between \(R_{n, \omega}\) and $H_{\L_n, \mathrm{per}}(\omega_{\L_n})$ is given in the next lemma. 

    \begin{lemma}\label[lemma]{le_energy_density_of_stationary_empirical_field}
            We have that
            \begin{align*}
                N_{\L_n}(\omega)
                = \abs{\L_n} R_{n, \omega}[N_C]\qquad\text{ and }\qquad
                H_{\L_n, \mathrm{per}}(\omega)
                = \abs{\L_n} H(R_{n, \omega}).
            \end{align*}
        \end{lemma}
        \begin{proof}
            The first equality is included in ~\cite[Remark 2.3, (1)]{GZ93}. Further, by~\cite[Remark 2.3, (3)]{GZ93}, we have that
            \begin{align*}
                (R_{n, \omega})^o
                = \abs{\L_n}^{-1} \sum_{(x, R_x) \in \omega_{\Lambda_n}} \delta_{(R_x, \theta_x \omega^{(n)})}
            \end{align*} is the Palm measure of \(R_{n, \omega}\).
            Therefore, \Cref{prop_existence_of_energy_density} implies
            \begin{align*}
                \abs{\L_n} H(R_{n, \omega})
                = \abs{\L_n} (R_{n, \omega})^o(f_H)
                = H_{\L_n, \mathrm{per}}(\omega),
            \end{align*} at least if the periodized configuration \(\omega^{(n)}\) respects the hardcore condition. But otherwise both left- and right-hand side are \(\infty\).
        \end{proof}

        \begin{proof}[Proof of \Cref{prop_pressure_per}]
            Let us first note that
            \begin{align*}
                I^z(P)
                = I(P) + (z-1) - \log(z)P(N_C)
            \end{align*} with \(I(P) := I^1(P)\) and second one has that
            \begin{align*}
                \frac{\log(\pi^z_{\Lambda_n}\left[\e^{-\beta H_{\Lambda_n, \mathrm{per}}} \right])}{\vert \Lambda_n \vert}
                &= \frac{\log(\pi_{\Lambda_n}\left[\e^{-\beta H_{\Lambda_n, \mathrm{per}}} e^{\abs{\L_n}(1-z)}z^{-N_{\L_n}} \right])}{\vert \Lambda_n \vert} \\
                &= \frac{\log(\pi_{\Lambda_n}\left[\e^{- (\beta H_{\Lambda_n, \mathrm{per}} - \log(z)N_{\L_n})} \right])}{\vert \Lambda_n \vert} + (1-z).
            \end{align*}
            Hence, the upper bound
            \begin{align*}
                \limsup_{n \to \infty} \frac{\log(\pi^z_{\Lambda_n}\left[\e^{-\beta H_{\Lambda_n, \mathrm{per}}} \right])}{\vert \Lambda_n \vert}
                &= \limsup_{n \to \infty} \frac{\log(\pi_{\Lambda_n}\left[\e^{- (\beta H_{\Lambda_n, \mathrm{per}} - \log(z)N_{\L_n})} \right])}{\vert \Lambda_n \vert} + (1-z)\\
                &\leq -\inf_{P \in \mathcal{P}_{\Theta}} \left[I(P) + \beta H(P) -\log(z)P(N_C)\right] + (1-z) \\
                &= -\inf_{P \in \mathcal{P}_{\Theta}} \left[I^z(P) + \beta H(P)\right]
            \end{align*}
            follows directly from~\cite[Theorem 3.1]{GZ93}. The assumptions of the cited theorem are met since \(H\) is trivially stable by being positive and lower-semicontinuous by \Cref{thm_lsc_of_energy_density} and we have \(\big\vert \abs{\L_n}^{-1} \bigl( \beta H_{\L_n, \mathrm{per}}(\omega) - \log(z) N_{\L_n} \bigr) - \bigl( \beta H(R_{n, \omega}) - \log(z) R_{n, \omega}[N_C] \bigr) \big\vert = 0\) by \Cref{le_energy_density_of_stationary_empirical_field}.
            
            The lower bound
            \begin{align*}
                \liminf_{n \to \infty} \frac{\log(\pi^z_{\Lambda_n}\left[\e^{-\beta H_{\Lambda_n, \mathrm{per}}} \right])}{\vert \Lambda_n \vert}
                \geq -\inf_{P \in \mathcal{P}_{\Theta}} \left[I^z(P) + \beta H(P)\right]
            \end{align*} does not directly follow from the mentioned theorem of Georgii--Zessin, but it follows from a closer inspection of its proof, together with the fact that \(H\) is continuous on \(\{H < \infty\}\) (see \Cref{thm_lsc_of_energy_density}) and the ergodic approximation given in \Cref{le_ergodic_approximation} below.
        \end{proof}

        \begin{lemma}\label[lemma]{le_ergodic_approximation}
            Let \(P \in \mathcal{P}_{\Theta}\) with \(I^z(P) + H(P) <\infty\).
            Then, there is a sequence \((P^{(n)})_{n \in \mathbb{N}_0}\) of ergodic probability measures with \(P^{(n)} \xrightarrow[n \uparrow \infty]{} P\) in the \(\tau_{\mathcal{L}}\)-topology and
            \begin{align*}
                \limsup_{n \to \infty} I^z(P^{(n)}) + H(P^{(n)}) \leq I^z(P) + H(P) 
            \end{align*}
        \end{lemma}
        \begin{proof}
            Consider the events \(B_n = \{\omega \in \Omega_{\L_{n}} \,\vert\, \forall \x \in \omega_{\L_n} \colon R_x < n^{1 - \epsilon}\}\), \(n \in \mathbb{N}\), for some fixed \(\epsilon \in (0, 1)\). Let \(m_n := n + n^{1-\epsilon}\), \(n \in \mathbb{N}\).
            Define \(\widehat{P}^{(n)}\) to be the probability measure on the configuration space for which the restrictions to \(\L_{m_n} + 2m_n i\), \(i \in \mathbb{Z}^d\), are iid with a distribution which does not place particles in \((\L_{m_n}\setminus \L_n) + 2m_n i\) and equals \(P(\cdot \,\vert\, B_n) \circ \theta_i^{-1}\) on \(\L_{n} + 2m_n i\).
            Notice that \(\widehat{P}^{(n)}\) respects the hardcore condition by construction (by restriction to \(B_n\)) and it inherits from $P$ the fact that it has a finite energy density.
            Now let
            \begin{align*}
                P^{(n)} := \frac{1}{\abs{\L_{m_n}}} \int_{\L_{m_n}} \, (\widehat{P}^{(n)} \circ \theta_x^{-1}) \, \mathrm{d}x.
            \end{align*}
            It is somewhat standard to show that each \(P^{(n)}\) is ergodic, that \(\limsup_{n \to \infty} I^z(P) \leq I^z(P)\), and \(P^{(n)} \xrightarrow[n \uparrow \infty]{} P\) in the \(\tau_{\mathcal{L}}\)-topology, see for example~\cite[Lemma 5.1]{Georgii_1994} for a similar construction and details in this analogous case.
            The only two additional ingredients that we need are $P(B_n) \xrightarrow[n \uparrow \infty]{} 1$ and that $\sup_{x \in \mathbb{R}^d} P_{\Lambda_n}\big[N_{x + \Delta}^{\psi} \mathbbm{1}_{B_n^c} \big] \xrightarrow[n \uparrow \infty]{} 0$ for all bounded measurable \(\Delta \subseteq \mathbb{R}^d\).

            In \Cref{le_temperedness_estimate_poisson} we saw that
            \begin{align*}
                \pi^z_{\L_n}\left[B_n^c \right]
                \leq c \e^{-\abs{\L_n}^{1+\gamma}}
            \end{align*} for some \(\gamma = \gamma(\epsilon) > 0\) and all \(n\) large enough.

            Then we have, by inequality \Cref{eq_entropic_inequality},
            \begin{align*}
                P_{\Lambda_n}(B_n^c)
                &= \frac{P_{\Lambda_n}\big[\abs{\L_n}^{1+\gamma} \mathbbm{1}_{B_n^c}\big]}{\abs{\L_n}^{1+\gamma}}
                \leq \frac{\log\big(\pi^z_{\L_n}\big[\e^{\abs{\L_n}^{1+\gamma} \mathbbm{1}_{B_n^c}} \big] \big)}{\abs{\L_n}^{1+\gamma}} + \frac{\abs{\L_n}}{\abs{\L_n}^{1+\gamma}} I(P) \\
                &= \frac{\log\big(\e^{\abs{\L_n}^{1+\gamma}}\pi^z_{\L_n}(B_n^c) + \pi^z_{\L_n}(B_n) \big)}{\abs{\L_n}^{1+\gamma}} + \frac{I^z(P)}{\abs{\L_n}^{\gamma}} \\
                &\leq \frac{\log\left(c +1\right)}{\abs{\L_n}^{1+\gamma}} + \frac{I^z(P)}{\abs{\L_n}^{\gamma}}
                \xrightarrow[n \uparrow \infty]{} 0.
            \end{align*}
            With the same entropy inequality, we see that \((N_{x+\Delta}^{\psi})_{x \in \mathbb{R}^d}\) is uniformly integrable w.r.t.~\(P\). Indeed, this follows with \(N_{x+\Delta}^{\psi}(\omega) = N_{x+\Delta}(\omega) + \sum_{y \in \omega_{\Delta + x}} R_y^{d}\) and
            \begin{align*}
P\left[\left(N_{x+\Delta}^{\psi}\right)^{1+\frac{d}{\delta}} \right]
                &= P\left[\left(N_{x+\Delta}^{\psi}\right)^{1+\frac{d}{\delta}} \right]
                \\
                &= P\left[\left(N_{x+\Delta}^{\psi}\right)^{1+\frac{d}{\delta}} \mathbbm{1}_{N_{x+\Delta} \leq 1}\right] + P\left[\left(N_{x+\Delta}^{\psi}\right)^{1+\frac{d}{\delta}} \mathbbm{1}_{N_{x+\Delta} > 1}\right]\\
                &\lesssim \vert x + \Delta \vert + P\Big[\sum_{y \in \omega_{\Delta + x}} R_y^{d + \delta}\Big] \\
                &\leq \vert \Delta \vert (1 + I^z(P)) + \log\Big(\pi^z\Big[\e^{\sum_{y \in \omega_{\Delta + x}} R_y^{d + \delta}} \Big] \Big)\\
                &= \vert \Delta \vert (I^z(P) + (1 - z) + z \mathfrak{R}(\e^{R^{d+\delta}}))
                < \infty,
            \end{align*} where we used the hardcore condition in the first inequality, and the analogous and well-known bound
            \begin{equation*}
                \sup_{x \in \mathbb{R}^d}P[N_{x+\Delta} \cdot \log_+(N_{x+\Delta})] \leq \abs{\Delta}I^z(P) + \log(\pi^z[\e^{N_{\Delta} \log_+(N_{\Delta})}]) < \infty.
            \end{equation*}
            Hence, for all bounded measurable \(\Delta \subseteq \mathbb{R}^d\),
            \begin{equation*}
                \sup_{x \in \mathbb{R}^d} P_{\Lambda_n}\big[N_{x + \Delta}^{\psi} \mathbbm{1}_{B_n^c} \big] \leq \sup_{x \in \mathbb{R}^d} P\big[N_{x + \Delta}^{\psi} \mathbbm{1}_{B_n^c} \big]  \xrightarrow[n \uparrow \infty]{} 0.  
            \end{equation*}
            Finally, as \(H\) is continuous on the set \(\{H < \infty\}\), and \(H(P^{(n)}) < \infty\), we have that \(\limsup_{n \to \infty} I^z(P^{(n)}) + H(P^{(n)}) \leq I^z(P) + H(P)\).
        \end{proof}

    \subsection{Free boundary conditions}
    Let us write $H_{\L,\emptyset}(\omega)=H_{\L}(\omega_\L)$ for the energy with free boundary conditions. The aim of this section is the following extension of \Cref{prop_pressure_per} to free boundary  conditions. 
        \begin{proposition}\label[proposition]{prop_pressure_free}
                We have
                \begin{align*}
                    \lim_{n \to \infty} \frac{\log(\pi^z_{\Lambda_n}\left[\e^{-\beta H_{\Lambda_n, \emptyset}} \right])}{\vert \Lambda_n \vert}
                    = -\inf_{P \in \mathcal{P}_{\Theta}} \left[I^z(P) + \beta H(P)\right].
                \end{align*}
            \end{proposition}

    As already explained above, this is a direct consequence of the comparison estimate in \Cref{le_comparison_pressure_per_and_free} which we will now prove. 
    
    \begin{proof}[Proof of \Cref{le_comparison_pressure_per_and_free}]
                Remember that
                \begin{align}\label{eq:recall_temperedness_estimate_poisson}
                    \pi^z\bigl( \exists \x \in \omega_{\L_K} \colon R_x \geq g(K) \bigr)
                    = \pi^z((\Omega^*_{K, K})^c)
                    \lesssim \e^{-\abs{\L_K}^{1+\gamma}}
                \end{align} by \Cref{le_temperedness_estimate_poisson}, which we will use in the following to establish limits.
                Consider the events \(A_n = \{\omega \in \Omega_{\L_{n+g(n)}} \,\vert\, \omega_{\L_{n+g(n)} \setminus \L_{n}} = \emptyset\}\) and \(B_n = \{\omega \in \Omega_{\L_{n+g(n)}} \,\vert\, \forall \x \in \omega_{\L_n} \colon R_x < g(n)\}\). Then, on \(A_n \cap B_n\), we have \(H_{\L_{n+g(n)},\mathrm{per}} =  H_{\L_{n+g(n)}, \emptyset} = H_{\L_{n}, \emptyset}\) and therefore
                \begin{align*}
                    \pi^z_{\L_{n+g(n)}}&\Big[\e^{-\beta H_{\L_{n+g(n)},\mathrm{per}}} \Big]
                    \geq \pi^z_{\L_{n+g(n)}}\Big[\e^{-\beta H_{\L_{n+g(n)}, \mathrm{per}}} \mathbbm{1}_{A_n}\mathbbm{1}_{B_n}\Big]\\
                    &=
                    \pi^z_{\L_{n+g(n)}}\Big[\e^{-\beta H_{\L_{n+g(n)}, \emptyset}} \mathbbm{1}_{A_n}\mathbbm{1}_{B_n}\Big] \\
                    &=
                    \pi^z_{\L_{n+g(n)}}\Big[\e^{-\beta H_{\L_{n+g(n)}, \emptyset}} \mathbbm{1}_{A_n}\Big] - \pi^z_{\L_{n+g(n)}}\Big[\e^{-\beta H_{\L_{n+g(n)}, \emptyset}} \mathbbm{1}_{A_n}\mathbbm{1}_{{B_n}^c}\Big]
                    \\
                    &=
                    \e^{z\abs{\L_{n + g(n)}} - z\abs{\L_n}}\pi^z_{\L_{n}}\Big[\e^{-\beta H_{\L_{n}, \emptyset}} \Big]
                    - \pi^z_{\L_{n+g(n)}}\Big[\e^{-\beta H_{\L_{n+g(n)}, \emptyset}} \mathbbm{1}_{A_n}\mathbbm{1}_{{B_n}^c}\Big] 
                \end{align*}
                It follows, using the estimate \Cref{eq:recall_temperedness_estimate_poisson}, that
                \begin{align*}
                    \limsup_{n \to \infty} \frac{\log\big(\pi^z_{\L_n}\big[\e^{-\beta H_{\L_n, \mathrm{per}}} \big])}{\vert \L_n \vert}
                    &\geq
                    \limsup_{n \to \infty} \frac{\log\Big(\pi^z_{\L_{n+g(n)}}\Big[\e^{-\beta H_{\L_{n+g(n)}, \mathrm{per}}} \Big]\Big)}{\vert \L_{n+g(n)} \vert}\\
                    &\geq 
                    \limsup_{n \to \infty} \frac{\log(\pi^z_{\L_n}\left[\e^{-\beta H_{\L_n, \emptyset}} \right])}{\vert \L_n \vert}.
                \end{align*}
                Finally, considering the event \(B'_n = \{\omega \in \Omega_{\L_{n}} \,\vert\, \forall \x \in \omega_{\L_n} \colon R_x < g(n)\}\), we have
                \begin{align*}
                    \e^{-\beta H_{\L_n, \emptyset}} \mathbbm{1}_{B'_n}&\geq \Big(\e^{\beta c\big(\abs{\L_n} - \abs{\L_{n - g(n)}}\big)} \e^{-\beta H_{\L_n, \mathrm{per}}}\mathbbm{1}_{H_{\L_n, \mathrm{per}} < \infty} + 0 \cdot \mathbbm{1}_{H_{\L_n, \mathrm{per}} = \infty}\Big) \mathbbm{1}_{B'_n}\\
                    &= \e^{\beta c\big(\abs{\L_n} - \abs{\L_{n - g(n)}}\big)} \e^{-\beta H_{\L_n, \mathrm{per}}} \mathbbm{1}_{B'_n}
                \end{align*} for some \(c > 0\), which implies
                \begin{align*}
                    \pi^z_{\L_n}\Big[\e^{-\beta H_{\L_n, \emptyset}}\Big]
                    \geq \e^{\beta c\big(\abs{\L_n} - \abs{\L_{n - g(n)}}\big)}&\pi^z_{\L_n}\left[\e^{-\beta H_{\L_n, \mathrm{per}}} \right] + \pi^z_{\L_n}\Big[\e^{-\beta H_{\L_n, \emptyset}} \mathbbm{1}_{{B'_n}^c}\Big]\\
                    &- \e^{\beta c\big(\abs{\L_n} - \abs{\L_n - g(n)}\big)}\pi^z_{\L_n}\Big[\e^{-\beta H_{\L_n, \mathrm{per}}}  \mathbbm{1}_{{B'_n}^c}\Big], 
                \end{align*} 
                and thus, using again the estimate \Cref{eq:recall_temperedness_estimate_poisson},
                \begin{align*}
                    \liminf_{n \to \infty} \frac{\log\big(\pi^z_{\L_n}\left[\e^{-\beta H_{\L_n, \emptyset}} \right]\big)}{\vert \L_n \vert}
                    \geq 
                    \liminf_{n \to \infty} \frac{\log\big(\pi^z_{\L_n}\left[\e^{-\beta H_{\L_n, \mathrm{per}}} \right]\big)}{\vert \L_n \vert}
                \end{align*} as desired.
            \end{proof}

    \subsection{Tempered boundary conditions sampled from a finite-energy measure}
    As a last step, we now show how to extend the existence of the pressure to boundary conditions sampled from a tempered measure $P$ with finite energy density. These configurations will then of course always respect the hardcore condition.
        
        \begin{proof}[Proof of \Cref{prop_pressure_with_boundary_conditions}]
            The result will follow by comparing with the case of free boundary conditions, \Cref{prop_pressure_free}.
            Analogously to the proof of \Cref{le_comparison_pressure_per_and_free}, consider the event  \[A_n = \{\omega \in \Omega_{\L_n} \,\vert\, \omega_{\L_{n} \setminus \L_{n - g(n) - g(n+1)}} = \emptyset\},\]
            here with a slightly bigger \textit{buffer zone}, and \[B_n = \{\omega \in \Omega_{\L_{n}} \,\vert\, \forall \x \in \omega_{\L_n} \colon R_x < g(n)\}.\]

            The temperedness assumption on \(P\) means that we can assume that there is some \(N \in \mathbb{N}\) (depending on \(\zeta\)), such that \(R_x \leq g(n)\) for all \(\x \in \zeta_{\Lambda_{n}}\) and \(n \geq N\).

            Then \(H_{\L_n, \zeta} = H_{\L_n, \emptyset}\) on \(A_n \cap B_n\) and therefore
            \begin{align*}
                    &\pi^z_{\L_n}\left[\e^{-\beta H_{\L_n, \zeta}} \right]
                    \geq \pi^z_{\L_n}\left[\e^{-\beta H_{\L_n, \zeta}} \mathbbm{1}_{A_n}\mathbbm{1}_{B_n} \right]
                    =
                    \pi^z_{\L_n}\left[\e^{-\beta H_{\L_{n - g(n)-g(n+1)}, \emptyset}} \mathbbm{1}_{A_n} \mathbbm{1}_{B_n} \right] \\
                    &=
                    \pi^z_{\L_n}\left[\e^{-\beta H_{\L_{n - g(n)-g(n+1)}, \emptyset}} \mathbbm{1}_{A_n} \right] 
                    - \pi^z_{\L_n}\left[\e^{-\beta H_{\L_{n - g(n)-g(n+1)}, \emptyset}} \mathbbm{1}_{A_n} \mathbbm{1}_{B_n^c} \right] 
                    \\
                    &=
                    \e^{z\abs{\L_n} - z\abs{\L_{n - g(n)-g(n+1)}}}\pi^z_{\L_{n - g(n)}}\left[\e^{-\beta H_{\L_{n - g(n)}, \emptyset}} \right]
                    - \pi^z_{\L_n}\left[\e^{-\beta H_{\L_{n - g(n)-g(n+1)}, \emptyset}} \mathbbm{1}_{A_n} \mathbbm{1}_{B_n^c} \right]. 
                \end{align*}
                It follows that
                \begin{align*}
                    \liminf_{n \to \infty} \frac{\log\big(\pi^z_{\L_n}\big[\e^{-\beta H_{\L_n, \zeta}} \big]\big)}{\vert \L_n \vert}
                    &\geq 
                    \liminf_{n \to \infty} \frac{\log\Big(\pi^z_{\L_{n-g(n)-g(n+1)}}\Big[\e^{-\beta H_{\L_{n-g(n)-g(n+1)}, \emptyset}} \Big]\Big)}{\vert \L_{n-g(n)} \vert}\\
                    &\geq 
                    \liminf_{n \to \infty} \frac{\log\big(\pi^z_{\L_n}\big[\e^{-\beta H_{\L_n, \emptyset}} \big]\big)}{\vert \L_n \vert}.
                \end{align*}
                Moreover, for some \(c > 0\),
                \begin{align*}
                    \e^{-H_{\L_n, \emptyset}} \mathbbm{1}_{B_n}
                    &\geq \left(\e^{c\big(\abs{\L_n} - \abs{\L_{n - g(n) - g(n+1)}}\big)} \e^{-H_{\L_n, \zeta}}\mathbbm{1}_{H_{\L_n, \zeta} < \infty} + 0 \cdot \mathbbm{1}_{H_{\L_n, \zeta} = \infty}\right) \mathbbm{1}_{B_n}\\
                    &= \e^{c\big(\abs{\L_n} - \abs{\L_{n - g(n) - g(n+1)}}\big)} \e^{-H_{\L_n, \zeta}} \mathbbm{1}_{B_n},
                \end{align*}
                which implies that
                \begin{align*}
                    \pi^z_{\L_n}\Big[\e^{-\beta H_{\L_n, \emptyset}}\Big]
                    &\geq \e^{\beta c\big(\abs{\L_n} - \abs{\L_{n - g(n) - g(n+1)}}\big)}\pi^z_{\L_n}\left[\e^{-\beta H_{\L_n, \zeta}} \right] + \pi^z_{\L_n}\Big[\e^{-\beta H_{\L_n, \emptyset}} \mathbbm{1}_{B_n^c}\Big] \\
                    &\phantom{\pi^z_{\L_n}\Big[\e^{-\beta H_{\L_n, \emptyset}}\Big]}
                    - \e^{\beta c\big(\abs{\L_n} - \abs{\L_{n - g(n) - g(n+1)}}\big)}\pi^z_{\L_n}\Big[\e^{-\beta H_{\L_n, \zeta}}  \mathbbm{1}_{B_n^c}\Big].
                \end{align*} 
                We then conclude that
                \begin{align*}
                    \liminf_{n \to \infty} \frac{\log\big(\pi^z_{\L_n}\big[\e^{-\beta H_{\L_n, \emptyset}} \big]\big)}{\vert \L_n \vert}
                    = \limsup_{n \to \infty} \frac{\log\big(\pi^z_{\L_n}\big[\e^{-\beta H_{\L_n, \emptyset}} \big]\big)}{\vert \L_n \vert}
                    \geq \limsup_{n \to \infty} \frac{\log\big(\pi^z_{\L_n}\big[\e^{-\beta H_{\L_n, \zeta}}\big]\big)}{\vert \L_n \vert},
                \end{align*}
                which completes the proof.
        \end{proof}

\section{Proof of the variational principle}\label{sec_variational_principle}
In this section, we put together the arguments needed to prove the variational principle, as stated in \Cref{thm_variational_principle}.

    \subsection{Existence and non-triviality of minimizers}
    \begin{proof}[Proof of \Cref{le_free_energy_attains_finite_infimum}]
        First note that \(P \mapsto I^z(P) + \beta H(P)\) takes finite values, since for \(P = \delta_{\emptyset}\) we see \(I(P) + \beta H(P) = z + \beta 0 = z < \infty\).
        We now argue that \(\{I^z + \beta H \leq c\}\) is compact (in the \(\tau_\mathcal{L}\)-topology) for all \(c \in \mathbb{R}\), which is enough to see that an infimum of \(P \mapsto I^z(P) + \beta H(P)\) has to be attained.
       For this, note that the level sets \(\{I^z \leq c\}\) are compact and sequentially compact for all \(c \in \mathbb{R}\), see \Cref{prop_specific_entropy_basics}. 
        Moreover, \(\{H \leq c\}\) is closed and \(H\) is continuous on this set for all \(c \in \mathbb{R}\), see the proof of \Cref{thm_lsc_of_energy_density} in \Cref{sec_energy_density}.
        As \(H\) is also positive, this already implies that \(\{I^z + \beta H \leq c\}\) is a closed subset of the compact set \(\{I^z \leq c\}\) and therefore also compact for all \(c \in \mathbb{R}\).
    \end{proof}

  \subsection{Gibbs measures are minimizers}\label{sec_gibbs_are_minimizers}
    We now show that infinite-volume translation-invariant Gibbs measures are minimizers of the variational formula. For this, we first show that the energy density does not depend on the choice of boundary conditions. This will allow us to compare restrictions of infinite-volume Gibbs measures to the finite-volume Gibbs measure with free boundary conditions.  
    \begin{lemma}\label[lemma]{le:energy_density_with_special_boundary_condition}
        Let \(P \in \mathcal{P}_{\Theta}\) with \(P\big(\sum_{\{\x, \y \} \subseteq \omega} \phi_{\mathrm{hc}}(\x, \y)\big) < \infty\) and \(P\big(\sum_{\x \in\omega_{C}} R_x^d\big) < \infty\).
        Then we have
        \begin{align*}
                \lim_{n \uparrow \infty} \frac{P\left(\left\vert H_{\L_n} - H_{\L_n, \emptyset}\right\vert\right)}{\vert \L_n \vert} 
                = 0.
            \end{align*}
    \end{lemma}
    \begin{proof}
        We have
        \begin{align*}
            \left\vert H_{\L_n}(\omega) - H_{\L_n, \emptyset}(\omega\right\vert
            &= H_{\L_n, \emptyset}(\omega) - H_{\L_n}(\omega)\\
            &= -\sum_{\x \in \omega_{\L_n}} \Big( \sum_{k \geq 2} \tfrac{1}{k} \sum_{\substack{\{\y_1, \dots, \y_{k-1}\} \subseteq \omega \setminus \x,\\ \exists j \in \{1, \dots, k-1\} \colon \y_j \in \omega_{\L_n^c}}} \phi_k(\x, \y_1, \dots, \y_{k-1}) \Big)
        \end{align*}
        and an analogous argument as in \Cref{le:representation_of_energy_in_finite_volume} and \Cref{le:representation_of_energy_in_infinite_volume} shows
        \begin{align*}
            &\frac{P(\vert H_{\L_n} - H_{\L_n, \emptyset}\vert)}{\vert \L_n \vert} \\
            &= - P^o\Big(\sum_{k \geq 2} \tfrac{1}{k} \sum_{\{\y_1, \dots, \y_{k-1}\} \subseteq \omega\setminus\oo} \phi_k(\oo, \y_1, \dots, \y_{k-1}) \frac{\abs{\L_n \cap (\L_n^c - y_1)}}{\abs{\L_n}} \Big)
            \xrightarrow[n \uparrow \infty]{} 0,
        \end{align*}
        as desired.
    \end{proof}

        This now enables us to show the first direction of the Gibbs variational principle. 
        \begin{proof}[Proof of \Cref{thm_gibbs_are_minimizers}]    
            We have
            \begin{align*}
                I(P_{\L_n} \,\vert\, G_{\L_n, z, \beta})
                &= P_{\L_n}\left[\log\left(\frac{\mathrm{d}P_{\L_n}}{\mathrm{d}G_{\L_n, z, \beta}}\right)\right]
                = P_{\L_n}\left[\log\left(Z_{\L_n, z, \beta} \e^{\beta H_{\L_n}, \emptyset}\frac{\mathrm{d}P_{\L_n}}{\mathrm{d}\pi^z_{\L_n}} \right)\right]\\
                &= \log(Z_{\L_n, z, \beta}) + \beta P_{\L_n}\left[H_{\L_n, \emptyset} \right] +  I(P_{\L_n} \,\vert\, \pi^z_{\L_n})
            \end{align*}
            and therefore
            \begin{equation}\label{eq_free_energy_difference}
                \begin{split}
                    \delta F(P)
                    &:= \lim_{n \to \infty} \frac{I(P_{\L_n} \,\vert\, G_{\L_n, z, \beta})}{\vert \L_n \vert}
                    = \lim_{n \to \infty} \frac{\log(Z_{\L_n, z, \beta})}{\abs{\L_n}} + I^z(P) + \beta H(P) \\
                    &= [I^z(P) + \beta H(P)] - \inf_{Q \in \mathcal{P}_{\Theta}} [I^z(Q) + \beta H(Q)]
                \end{split}
            \end{equation} by \Cref{prop_pressure_free}.
            We have to show that
            \begin{align*}
                \delta F(P) = 0.
            \end{align*}
            For this we upper bound $I(P_{\Lambda} \,\vert\, G_{\Lambda, z, \beta})$ for $\Lambda \subseteq \mathbb{R}^d$.
            From the DLR equations for $P$ it follows that
            \begin{align*}
                g_{\Lambda}(\omega)
                := \frac{\mathrm{d}P_{\Lambda}}{\mathrm{d}G_{\Lambda, z, \beta}}(\omega)
                = P \left[\frac{\mathrm{d}G_{\Lambda, z, \beta, \bullet}(\omega)}{\mathrm{d}G_{\Lambda, z, \beta}(\omega)} \right].
            \end{align*}
            Herein we have
            \begin{align*}
                \frac{\mathrm{d}G_{\Lambda, z, \beta, \zeta}(\omega)}{\mathrm{d}G_{\Lambda, z, \beta}(\omega)}
                = \frac{Z_{\Lambda, z, \beta}}{Z_{\Lambda, z, \beta, \zeta}} \, \, \e^{\beta H_{\Lambda, \emptyset}(\omega) - \beta H_{\L, \zeta}(\omega)} .
            \end{align*}
            With the densities $g_{\Lambda}$ calculated above it is 
            \begin{align*}
                I(P_{\Lambda} \,\vert\, G_{\Lambda, z, \beta})
                = G_{\Lambda, z, \beta} \left[g_{\Lambda} \log(g_{\Lambda}) \right].
            \end{align*}
            Since the function $x \mapsto x \log(x)$ is convex, we can apply Jensen's inequality to obtain
            \begin{align*}
                g_{\Lambda}(\omega) \log(g_{\Lambda}(\omega))
                &= P \left[\frac{\mathrm{d}G_{\Lambda, z, \beta, \bullet}(\omega)}{\mathrm{d}G_{\Lambda, z, \beta}(\omega)} \right] \log\left( P \left[\frac{\mathrm{d}G_{\Lambda, z, \beta, \bullet}(\omega)}{\mathrm{d}G_{\Lambda, z, \beta}(\omega)} \right] \right) \\
                &\leq P\left[\frac{\mathrm{d}G_{\Lambda, z, \beta, \bullet}(\omega)}{\mathrm{d}G_{\Lambda, z, \beta}(\omega)} \log\left(\frac{\mathrm{d}G_{\Lambda, z, \beta, \bullet}(\omega)}{\mathrm{d}G_{\Lambda, z, \beta}(\omega)}\right) \right].
            \end{align*}
            By using that $P$ satisfies the  DLR equations we get 
            \begin{align*}
                I(P_{\Lambda} \,\vert\, G_{\Lambda, z, \beta})
                &= G_{\Lambda, z, \beta} \left[g_{\Lambda} \log(g_{\Lambda}) \right]\\
                &\leq G_{\Lambda, z, \beta}^{\omega} \left[ P^{\zeta}\left[\frac{\mathrm{d}G_{\Lambda, z, \beta, \zeta}(\omega)}{\mathrm{d}G_{\Lambda, z, \beta}(\omega)} \log\left(\frac{\mathrm{d}G_{\Lambda, z, \beta, \zeta}(\omega)}{\mathrm{d}G_{\Lambda, z, \beta}(\omega)}\right) \right] \right] \\
                &=  P^{\zeta}\left[G_{\Lambda, z, \beta}^{\omega}\left[\frac{\mathrm{d}G_{\Lambda, z, \beta, \zeta}(\omega)}{\mathrm{d}G_{\Lambda, z, \beta}(\omega)} \log\left(\frac{\mathrm{d}G_{\Lambda, z, \beta, \zeta}(\omega)}{\mathrm{d}G_{\Lambda, z, \beta}(\omega)}\right) \right] \right]\\
                &= P^{\zeta}\left[G_{\Lambda,z, \beta, \zeta}^{\omega}\left[\log\left(\frac{\mathrm{d}G_{\Lambda, z, \beta, \zeta}(\omega)}{\mathrm{d}G_{\Lambda, z, \beta}(\omega)}\right) \right] \right]\\
                &= P\left[\log\left(\frac{\mathrm{d}G_{\Lambda, z, \beta, \bullet}(\bullet)}{\mathrm{d}G_{\Lambda, z, \beta}(\bullet)}\right) \right] \\
                &= \beta P\left[H_{\L, \emptyset}(\bullet) - H_{\Lambda}(\bullet) \right] + \left[\log(Z_{\Lambda, z, \beta}) - P(\log(Z_{\Lambda, z, \beta, \bullet})) \right],
            \end{align*} where we used the notation \(Q_2^{\omega}\left[Q_1^{\zeta}\left[f(\omega, \zeta)\right]\right] = \int \, Q_2(\d \omega) \int \, Q_1(\d \zeta) f(\omega, \zeta) \) for readability reasons.
            Moreover, by Fatou's lemma
            \begin{align*}
                &\limsup_{n \to \infty} \frac{I(P_{\Lambda_n} \,\vert\, G_{\Lambda_n, z, \beta})}{\vert \Lambda_n \vert} \\
                &\leq \beta \limsup_{n \to \infty} \frac{P\left[H_{\L_n, \emptyset}(\bullet) - H_{\L_n}(\bullet) \right]}{\vert \Lambda_n \vert}  + \left[\limsup_{n \to \infty} \frac{\log(Z_{\Lambda_n, z, \beta})}{\vert \Lambda_n \vert} - \liminf_{n \to \infty} \frac{P[\log(Z_{\Lambda_n, z, \beta, \bullet})]}{\vert \Lambda_n \vert} \right]\\\\
                &\leq \beta \limsup_{n \to \infty} \frac{P\left[H_{\L_n, \emptyset}(\bullet) - H_{\L_n}(\bullet) \right]}{\vert \Lambda_n \vert}  + \left[\limsup_{n \to \infty} \frac{\log(Z_{\Lambda_n, z, \beta})}{\vert \Lambda_n \vert} - P\left[\liminf_{n \to \infty}\frac{\log(Z_{\Lambda_n, z, \beta, \bullet})}{\vert \Lambda_n \vert}\right] \right].
            \end{align*}
            But now \Cref{le:energy_density_with_special_boundary_condition} and \Cref{prop_pressure_with_boundary_conditions} imply
            \begin{align*}
                \limsup_{n \to \infty} \frac{I(P_{\Lambda_n} \,\vert\, G_{\Lambda_n, z, \beta})}{\vert \Lambda_n \vert} \leq 0
            \end{align*} 
            which of course implies,  by non-negativity of the relative entropy,   $$\delta F(P) = \lim_{n \to \infty} \frac{I(P_{\Lambda_n} \,\vert\, G_{\Lambda_n, z, \beta})}{\vert \Lambda_n \vert} = 0,$$
           which concludes the proof.
        \end{proof}

    \subsection{Minimizers are Gibbs measures}\label{sec_minimizers_are_gibbs}
    Now comes the harder part of the Gibbs variational principle.  We want to show that minimizers of the variational formula are Gibbs measures in the sense of \Cref{def_infinite_volume_Gibbs_measure}.
    
    Before we start with the proof of \Cref{thm_minimizers_are_gibbs} we prove two technical results which will allow us to estimate certain error terms. 
    The first one provides us with a tool to localize by cutting of boundary conditions outside of a sufficiently large volume. 
        
        \begin{lemma}\label[lemma]{le_total_variation_distance}
            Let \(P \in \mathcal{P}_\Theta\) with \(I^z(P) < \infty\).
            Then, for all bounded measurable \(L \subset \mathbb{R}^d\) and for all $\delta > 0$ there is a $K \in \mathbb{N}$ such that for all $k \geq K$ it holds
            \begin{align*}
                \int \, \nu(\mathrm{d}\zeta) \, \left\Vert G_{L,z,\beta, \zeta_{\Lambda_k}} - G_{L, z, \beta, \zeta} \right\Vert_{\text{TV}} < \delta
            \end{align*} for $\nu = P$ and $\nu = PG_{\Lambda}$ for all bounded measurable $\Lambda$.
        \end{lemma}
        \begin{proof}
            Consider for \(K \in \mathbb{N}\) the events \(\Omega_K^*\) from \eqref{sets_of_temperedness_conditions}
            and define 
            \begin{align*}
                B_K
                := \{\omega \in \Omega_{L} \,\vert\, \forall \x \in \omega_{L} \colon R_x \leq g(K)\}.
            \end{align*}
            Then we have \(H_{L, \zeta}(\omega) = H_{L, \zeta_{\Lambda_k}}(\omega)\) for \(K\) large enough, all \(k \geq K\) and \(\omega_{L} \in B_K, \zeta \in \Omega^*_{K}\).
            It follows that
            \begin{align*}
                \left\Vert G_{L, z, \beta, \zeta_{\Lambda_k}} - G_{L, z, \beta, \zeta} \right\Vert_{\text{TV}} 
                &\leq 
                2 \frac{\pi^z_{L}\left(\vert \e^{-\beta H_{L, \zeta}} - \e^{-\beta H_{L, \zeta_{\Lambda_k}}}\vert\right)}{\max\big\{Z_{L, z, \beta, \zeta}, Z_{L, z, \beta, \zeta_{\Lambda_k}}\big\}}
                \leq 2 \frac{\pi^z_{L}\left(2 \cdot \1_{B_K^c} \right)}{e^{-z\abs{L}}} \\
                &= 4 \e^{z\abs{L}} \,\pi^z_{L}\left(B_K^c\right),
            \end{align*} at least if \(\zeta \in \Omega^*_{K}\).
            Hence,
            \begin{align*}
                \int \, \nu(\mathrm{d}\zeta) \, \big\Vert G_{L,z, \beta, \zeta_{\Lambda_k}} - G_{L, z, \beta, \zeta} \big\Vert_{\text{TV}} 
                \leq 
                4 \e^{z \abs{L}} \,\pi^z_{L}\left(B_K^c\right) \nu(\Omega^*_{K}) + 2 \nu({\Omega^*_{K}}^c).
            \end{align*}
            For example the proof of \Cref{le_comparison_pressure_per_and_free} shows that \(\pi^z_{L}\left(B_K^c\right) \to 0\) and \(\nu({\Omega^*_{K}}^c) \to 0\), as \(K \to \infty\),  uniformly in \(\L \in \mathfrak{B}_b(\mathbb{R}^d)\) see \Cref{le_temperedness_estimate_for_nu}.
        \end{proof}
        Now note that, if $0 =\delta F(P) = \lim_{n \to \infty} \frac{1}{\abs{\Lambda_n}}I(P_{\Lambda_n} \,\vert\, G_{\Lambda_n, z, \beta})$,
            then, by the translation invariance of $P$, it follows that $P_\Lambda \ll G_{\Lambda, z, \beta}$ with some density $g_\Lambda$ for large enough cubes (not necessarily centered at the origin) $\Lambda$.
            For every subset $\Delta \subset \Lambda$, we consider the restriction $$G_{\Lambda, \Delta} = (G_{\Lambda, z, \beta})_{\Delta}.$$
            From the DLR equations for Gibbs measures in finite volumes, we have $P_\Delta \ll G_{\Lambda, \Delta}$, with density
            \begin{align*}
                g_{\Lambda, \Delta}(\zeta_\Delta) := \int \, G_{\Lambda \setminus \Delta, z, \beta, \zeta_\Delta}(\mathrm{d}\omega_{\Lambda\setminus\Delta}) \, g_{\Lambda}(\omega_{\Lambda\setminus\Delta} \zeta_{\Delta}).
            \end{align*}
            We now estimate the expected cost of replacing $g_{\Lambda, \Delta}$ by $g_{\Lambda, \Delta\setminus L}$. 
            
        \begin{lemma}[Consequence of $\delta F(P) = 0$ for involved densitites]\label[lemma]{lem_consequence_of_deltaF_0}
            Consider the situation of (the proof of) \Cref{thm_minimizers_are_gibbs}, in particular a $P$ with $\delta F(P) = 0$.
            Then, for every every $k \in \mathbb{N}$ with \(L \subseteq \Lambda_k\) and $\delta > 0$, there exist sets $\Delta, \Lambda \subset \mathbb{R}^d$ with $\Lambda_k \subset \Delta \subset \Lambda$, such that 
            \begin{align*}
                G_{\Lambda, z, \beta}[\vert g_{\Lambda, \Delta} - g_{\Lambda, \Delta\setminus L} \vert] < \delta.
            \end{align*}
        \end{lemma}
        \begin{proof}
            Let us denote \(I(P_{\Delta} \,\vert\, (G_{\Lambda, z, \beta})_{\Delta})\) by \(I_{\Delta}(P \,\vert\, G_{\Lambda, z, \beta})\) for bounded measurable sets \(\Delta \subseteq \Lambda \subset \R^d\) in this proof.
            According to \Cref{lem_comparison_of_densities_and_relative_entropies} the claim follows if we can show that for all $\delta > 0$  there are sets $\Delta, \Lambda \subset \mathbb{R}^d$ with $\Lambda_k \subset \Delta \subset \Lambda$, such that \begin{align*}
                I_{\Delta}(P \,\vert\, G_{\Lambda, z, \beta}) - I_{\Delta\setminus L}(P \,\vert\, G_{\Lambda, z, \beta})
                < \delta.
            \end{align*}
            The assumption $0 =\delta F(P) = \lim_{n \to \infty} \frac{I(P_{\Lambda_n} \,\vert\, G_{\Lambda_n, z, \beta})}{\vert \Lambda_n \vert}$ implies that there exists $n \geq k$ such that $\frac{I(P_{\Lambda_n} \,\vert\, G_{\Lambda_n, z, \beta})}{\vert \Lambda_n \vert} < \frac{\delta}{2^d \vert \Lambda_k \vert}$.
            Choose $m \in \mathbb{N}$ such that $k \leq \frac{n}{m} \leq 2k$ (e.g., $m = \lfloor \frac{n}{k} \rfloor$). We then have $m^d \vert \Lambda_k \vert \leq \vert \Lambda_n \vert \leq (2m)^d \vert \Lambda_k \vert$.
            Now choose $m^d$-many disjoint translates $\Lambda_k(l) = \Lambda_k + i(l) \subset \Lambda_n$ of $\Lambda_k$ in $\Lambda_n$, $i(l) \in \mathbb{Z}^d, 1 \leq l \leq m^d$, and denote $\Delta(l) = \bigcup_{i = 1}^{l} \Lambda_k(i)$.
            \medbreak
            
            We now justify why the claim follows by choosing \(\Lambda := \Lambda_n - i(l)\) and \(\Delta := \Delta(l) - i(l)\) for a specifically chosen \(l\).
            From the monotonicity of the relative entropy, we have 
            \begin{align*}
                &m^{-d} \sum_{l = 1}^{m^d} \left[I_{\Delta(l)}(P \,\vert\, G_{\Lambda_n, z, \beta}) - I_{\Delta(l) \setminus (L + i(l))}(P \,\vert\, G_{\Lambda_n, z, \beta}) \right] \\
                &\leq m^{-d} \sum_{l = 1}^{m^d} \left[I_{\Delta(l)}(P \,\vert\, G_{\Lambda_n, z, \beta}) - I_{\Delta(l) \setminus \Lambda_k(l)}(P \,\vert\, G_{\Lambda_n, z, \beta}) \right].
            \end{align*}
            Furthermore, using \(\Delta(l) = \Delta(l+1) \setminus \Lambda_k(l+1)\) for \(0 \leq l \leq m^{d} - 1\), 
            \begin{align*} 
                \sum_{l = 1}^{m^d} I_{\Delta(l)}(P \,\vert\, G_{\Lambda_n, z, \beta}) 
                &=  \sum_{l = 0}^{m^d - 1} I_{\Delta(l+1) \setminus \Lambda_k(l+1)}(P \,\vert\, G_{\Lambda_n, z, \beta}) + I_{\Delta(m^d)}(P \,\vert\, G_{\Lambda_n, z, \beta}) \\ 
                &= \sum_{l = 1}^{m^d} I_{\Delta(l) \setminus \Lambda_k(l)}(P \,\vert\, G_{\Lambda_n, z, \beta}) + I_{\Delta(m^d)}(P \,\vert\, G_{\Lambda_n, z, \beta}).
            \end{align*}
            This means (using the monotonicity of the entropy in the second inequality) 
            \begin{align*} 
                &m^{-d} \sum_{l = 1}^{m^d} \left[I_{\Delta(l)}(P \,\vert\, G_{\Lambda_n, z, \beta}) - I_{\Delta(l) \setminus (L + i(l))}(P \,\vert\, G_{\Lambda_n, z, \beta}) \right] \\
                &\leq m^{-d} I_{\Delta(m^d)}(P \,\vert\, G_{\Lambda_n, z, \beta}) 
                \leq m^{-d} I_{\Lambda_n}(P \,\vert\, G_{\Lambda_n, z, \beta}) \\
                &< \frac{\delta \vert\Lambda_n\vert}{(2m)^d \vert \Lambda_k \vert} 
                \leq \delta.
            \end{align*}
            Since the terms $I_{\Delta(l)}(P \,\vert\, G_{\Lambda_n, z, \beta}) - I_{\Delta(l) \setminus (L + i(l))}(P \,\vert\, G_{\Lambda_n, z, \beta})$ are on average smaller than $\delta$, there has to be at least one specific $l$ such that 
            \begin{equation*}
                I_{\Delta(l)}(P \,\vert\, G_{\Lambda_n, z, \beta}) - I_{\Delta(l) \setminus (L + i(l))}(P \,\vert\, G_{\Lambda_n, z, \beta}) < \delta.
            \end{equation*}
            Using translation invariance (i.e.,\(I_{\Delta + x}(P \,\vert\, G_{\Lambda + x, z, \beta}) = I_{\Delta}(P \,\vert\, G_{\Lambda, z, \beta}) \)), we have thus found suitable sets \(\Lambda := \Lambda_n - i(l)\) and \(\Delta := \Delta(l) - i(l)\).
        \end{proof}

        \begin{lemma}\label[lemma]{lem_comparison_of_densities_and_relative_entropies}
            Consider the situation of (the proof of) \Cref{lem_consequence_of_deltaF_0}.
            Then, for all $\delta > 0$, there exists a $\delta' > 0$ such that $G_{\Lambda, z, \beta}\left[\left\vert g_{\Lambda, \Delta} - g_{\Lambda, \Delta\setminus L} \right\vert\right] < \delta$ if $I_{\Delta}(P \,\vert\, G_{\Lambda, z, \beta}) - I_{\Delta\setminus L}(P \,\vert\, G_{\Lambda, z, \beta}) < \delta'$.
        \end{lemma}

        \begin{proof}
            Consider the function $\psi(x) = 1 - x + x \log(x)$.
            Then, for any $\epsilon > 0$, there exists $r_{\epsilon} > 0$ such that
            \begin{equation*}
                \psi(x) \geq (\vert x -1 \vert - \epsilon/2)/r_{\epsilon}.
            \end{equation*}
            We compute
            \begin{align*} 
                &I_{\Delta}(P \,\vert\, G_{\Lambda, z, \beta}) - I_{\Delta\setminus L}(P \,\vert\, G_{\Lambda, z, \beta}) 
                = P[\log(g_{\Lambda, \Delta})] - P[\log(g_{\Lambda, \Delta\setminus L})] \\
                &= P\left[\log\left(\tfrac{g_{\Lambda, \Delta}}{g_{\Lambda, \Delta\setminus L}}\right)\right]
                = G_{\Lambda, z, \beta}\left[g_{\Lambda, \Delta} \log\left(\tfrac{g_{\Lambda, \Delta}}{g_{\Lambda, \Delta\setminus L}}\right)\right] \\
                &= 1 - 1 + G_{\Lambda, z, \beta}\left[g_{\Lambda, \Delta \setminus L} \left(\tfrac{g_{\Lambda, \Delta}}{g_{\Lambda, \Delta\setminus L}}\right)\log\left(\tfrac{g_{\Lambda, \Delta}}{g_{\Lambda, \Delta\setminus L}}\right)\right] \\
                &= G_{\Lambda, z, \beta}[g_{\Lambda, \Delta \setminus L}]  - G_{\Lambda, z, \beta}\left[g_{\Lambda, \Delta \setminus L} \left(\tfrac{g_{\Lambda, \Delta}}{g_{\Lambda, \Delta\setminus L}}\right)\right] + G_{\Lambda, z, \beta}\left[g_{\Lambda, \Delta \setminus L} \left(\tfrac{g_{\Lambda, \Delta}}{g_{\Lambda, \Delta\setminus L}}\right)\log\left(\tfrac{g_{\Lambda, \Delta}}{g_{\Lambda, \Delta\setminus L}}\right)\right] \\ 
                &= G_{\Lambda, z, \beta}\left[g_{\Lambda, \Delta\setminus L} \,\, \psi\left(\tfrac{g_{\Lambda, \Delta}}{g_{\Lambda, \Delta\setminus L}}\right)\right].
            \end{align*}
            From the above estimate for $\psi$, we get
            \begin{align*} 
                I_{\Delta}(P \,\vert\, G_{\Lambda, z, \beta}) - I_{\Delta\setminus L}(P \,\vert\, G_{\Lambda, z, \beta}) 
                &= G_{\Lambda, z, \beta}\left[g_{\Lambda, \Delta\setminus L} \,\, \psi\left(\tfrac{g_{\Lambda, \Delta}}{g_{\Lambda, \Delta\setminus L}}\right)\right] \\
                &\geq \frac{1}{r_\epsilon} \left( G_{\Lambda, z, \beta}\left[g_{\Lambda, \Delta\setminus L} \,\, \left\vert\tfrac{g_{\Lambda, \Delta}}{g_{\Lambda, \Delta\setminus L}} - 1\right\vert\right] - \tfrac{\epsilon}{2} \right) \\ 
                &= \frac{1}{r_\epsilon} \left( G_{\Lambda, z, \beta}\left[\left\vert g_{\Lambda, \Delta} - g_{\Lambda, \Delta\setminus L} \right\vert\right] - \tfrac{\epsilon}{2} \right),
            \end{align*}
            that is,
            \begin{align*}
                G_{\Lambda, z, \beta}\left[\left\vert g_{\Lambda, \Delta} - g_{\Lambda, \Delta\setminus L} \right\vert\right]
                \leq
                \tfrac{\epsilon}{2} + r_{\epsilon} \left(I_{\Delta}(P \,\vert\, G_{\Lambda, z, \beta}) - I_{\Delta\setminus L}(P \,\vert\, G_{\Lambda, z, \beta})\right).
            \end{align*}
            Choosing $\epsilon < \delta$ and $\delta' = \delta/(2 r_\epsilon)$ concludes the proof.
        \end{proof}

        With these technical tools at our hands, we are finally ready to provide the proof of \Cref{thm_minimizers_are_gibbs} and thereby finishing the proof of \Cref{thm_variational_principle}.

        \begin{proof}[Proof of \Cref{thm_minimizers_are_gibbs}]
            Let \(P \in \mathcal{P}_\Theta\) with \(I^z(P) + \beta H(P) = \inf_{Q \in \mathcal{P}_{\Theta}} \left[I^z(Q) + \beta H(Q)\right]\). As we have already seen in \Cref{eq_free_energy_difference}, this is equivalent to \(\delta F(P) = 0\). 
            
            Recall, that in order to show that $P$ is a Gibbs measure, we need to show that the DLR equation \Cref{eq:DLR} is satisfied for any bounded measurable $L \subset \R^d$ and any measurable $f\geq 0$, namely
            \begin{align*}
                \int \, P(\mathrm{d}\zeta) \, f(\zeta) 
                = \int \, P(\mathrm{d}\zeta) \, \left[\int \, G_{L, z, \beta, \zeta}(\mathrm{d}\omega) \, f(\omega_{L}\zeta_{L^c}) \right].
            \end{align*}
            
            Without loss of generality (using the functional monotone class theorem), we can assume $0 \leq f \leq 1$ and that $f$ is a local function, i.e., it only depends on the particles in a bounded region.
            We will show that, for any $\delta > 0$,
            \begin{equation}\label{eq:DLR_estimate}
                \left\vert \int \, P(\mathrm{d}\zeta) \, f(\zeta) - \int \, P(\mathrm{d}\zeta) \, \left[\int \, G_{L, z, \beta, \zeta}(\mathrm{d}\omega) \, f(\omega_{L}\zeta_{L^c}) \right] \right\vert 
                < 3 \delta.
            \end{equation}
            
            \medbreak
            
            First note that $0 =\delta F(P) = \lim_{n \to \infty} \frac{1}{\abs{\Lambda_n}}I(P_{\Lambda_n} \,\vert\, G_{\Lambda_n, z, \beta})$
            and the translation invariance of $P$, imply that $P_\Lambda \ll G_{\Lambda, z, \beta}$ with some density $g_\Lambda$ for sufficiently large cubes $\Lambda$.
            For every subset $\Delta \subset \Lambda$, we consider the restriction $$G_{\Lambda, \Delta} = (G_{\Lambda, z, \beta})_{\Delta}.$$
            From the DLR equations for Gibbs measures in finite volume, we have $P_\Delta \ll G_{\Lambda, \Delta}$, with density
            \begin{align*}
                g_{\Lambda, \Delta}(\zeta_\Delta) := \int \, G_{\Lambda \setminus \Delta, z, \beta, \zeta_\Delta}(\mathrm{d}\omega_{\Lambda\setminus\Delta}) \, g_{\Lambda}(\omega_{\Lambda\setminus\Delta} \zeta_{\Delta}).
            \end{align*}

            We are then able to choose $k$ so large that 
            \begin{enumerate}
                \item[(a)]\(L \subseteq \Lambda_k\),
                \item[(b)] $f$ just depends on particles inside $\Lambda_k$, 
                \item[(c)] and \begin{align*}
                \int \, \nu(\mathrm{d}\zeta) \, \big\Vert G_{L, z, \beta, \zeta_{\Lambda_k}} - G_{L, z, \beta, \zeta} \big\Vert_{\text{TV}} 
                < \delta
            \end{align*}
            for $\nu = P$ and $\nu = PG_{\Lambda}$ for all $\Lambda$.
            \end{enumerate}
            Here, $(a)$ and $(b)$ are clear and $(c)$ follows from \Cref{le_total_variation_distance}.

            Moreover, for each such $k$, we can apply \Cref{lem_consequence_of_deltaF_0} to choose $\Delta, \Lambda \subset \mathbb{R}^d$ such that $\Lambda_k \subset \Delta \subset \Lambda$, and
            \begin{align*}
                G_{\Lambda}\left[\left\vert g_{\Lambda, \Delta} - g_{\Lambda, \Delta\setminus L} \right\vert\right]
                < \delta.
            \end{align*}
            The rest of the proof is structured according to the following six steps:
            \begin{align*} 
                \int \, P(\mathrm{d}\zeta) \, f(\zeta)
                &\stackrel{(1)}{=} \int \, G_{\Lambda, z, \beta}(\mathrm{d}\zeta) \, g_{\Lambda, \Delta}(\zeta) \, f(\zeta) \\
                &\stackrel{(2)}{\approx} \int \, G_{\Lambda, z, \beta}(\mathrm{d}\zeta) \, g_{\Lambda, \Delta\setminus L}(\zeta) \, f(\zeta)\\
                &\stackrel{(3)}{=} \int \, G_{\Lambda, z, \beta}(\mathrm{d}\zeta) \, g_{\Lambda, \Delta\setminus L}(\zeta) \, \left[\int \, G_{L, z, \beta, \zeta}(\mathrm{d}\omega) \, f(\omega_{L} \zeta_{L^c}) \right]\\ 
                &\stackrel{(4)}{\approx} \int \, G_{\Lambda, z, \beta}(\mathrm{d}\zeta) \, g_{\Lambda, \Delta\setminus L}(\zeta) \, \left[\int \, G_{L, z, \beta, \zeta_{\Lambda_k}}(\mathrm{d}\omega) \, f(\omega_{L} \zeta_{\Lambda_{k}\setminus L}) \right] \\
                &\stackrel{(5)}{=} \int \, P(\mathrm{d}\zeta) \, \left[\int \, G_{L, z, \beta, \zeta_{\Lambda_k}}(\mathrm{d}\omega) \, f(\omega_{L} \zeta_{\Lambda_{k}\setminus L}) \right] \\ 
                &\stackrel{(6)}{\approx} \int \, P(\mathrm{d}\zeta) \, \left[\int \, G_{L, z, \beta, \zeta}(\mathrm{d}\omega) \, f(\omega_{L} \zeta_{\Lambda_{k}\setminus L}) \right],
            \end{align*}
            where we write $X \approx Y$ if $\vert X - Y \vert < \delta$.
            \medbreak
            
            \begin{description}[before={\renewcommand\makelabel[1]{\bfseries (##1)}},leftmargin=0pt,itemsep=10pt] 
                \item[1] This is true since
                \begin{align*} 
                    \int \, P(\mathrm{d}\zeta) \, f(\zeta)
                    = \int \, P(\mathrm{d}\zeta) \, f(\zeta_{\Delta})
                    = \int \, P_{\Delta}(\mathrm{d}\zeta) \, f(\zeta)
                    = \int \, G_{\Lambda, \Delta}(\mathrm{d}\zeta) \, g_{\Lambda, \Delta}(\zeta) \, f(\zeta) \\
                    = \int \, G_{\Lambda, z, \beta}(\mathrm{d}\zeta) \, g_{\Lambda, \Delta}(\zeta_{\Delta}) \, f(\zeta_{\Delta})
                    = \int \, G_{\Lambda, z, \beta}(\mathrm{d}\zeta) \, g_{\Lambda, \Delta}(\zeta) \, f(\zeta),
                \end{align*} 
                as $f$ only depends on particles inside $\Lambda_k \subset \Delta$ and $g_{\Lambda, \Delta}(\zeta) = \frac{\mathrm{d}P_{\Delta}(\zeta)}{\mathrm{d}G_{\Lambda, \Delta}(\zeta)}$ per definition.
                
                \item[2] This holds since, from the assumptions, we have $G_{\Lambda, z, \beta}[\vert g_{\Lambda, \Delta} - g_{\Lambda, \Delta\setminus L} \vert] < \delta$. Then, it immediately follows that
                \begin{align*}
                    &\left\vert \int \, G_{\Lambda, z, \beta}(\mathrm{d}\zeta) \, g_{\Lambda, \Delta}(\zeta) \, f(\zeta) - \int \, G_{\Lambda, z, \beta}(\mathrm{d}\zeta) \, g_{\Lambda, \Delta\setminus L
                    }(\zeta) \, f(\zeta)\right\vert \\
                    &\leq \int \, G_{\Lambda, z, \beta}(\mathrm{d}\zeta) \, \left\vert g_{\Lambda, \Delta}(\zeta) - g_{\Lambda, \Delta\setminus L}(\zeta) \right\vert \, \left\vert f(\zeta)\right\vert \\
                    &\leq G_{\Lambda, z, \beta}[\vert g_{\Lambda, \Delta} - g_{\Lambda, \Delta\setminus L} \vert] < \delta,
                \end{align*} 
                because we assumed $0 \leq f \leq 1$. 
                \item[3] This is true since,
                \begin{align*}
                    \int \, G_{\Lambda, z, \beta}(\mathrm{d}\zeta) \, g_{\Lambda, \Delta\setminus L}(\zeta) \, f(\zeta)
                    &= \int \, G_{\Lambda, z, \beta}(\mathrm{d}\zeta) \, \left[\int \, G_{L, z, \beta, \zeta}(\mathrm{d}\omega) \, g_{\Lambda, \Delta\setminus L}(\omega_{L} \zeta_{L^c}) \, f(\omega_{L} \zeta_{L^c}) \right] \\
                    &= \int \, G_{\Lambda, z, \beta}(\mathrm{d}\zeta) \, \left[\int \, G_{L, z, \beta, \zeta}(\mathrm{d}\omega) \, g_{\Lambda, \Delta\setminus L}(\zeta_{\Delta\setminus L}) \, f(\omega_{L} \zeta_{L^c}) \right] \\
                    &= \int \, G_{\Lambda, z, \beta}(\mathrm{d}\zeta) \, g_{\Lambda, \Delta\setminus L}(\zeta_{\Delta\setminus L}) \left[\int \, G_{L, z, \beta, \zeta}(\mathrm{d}\omega) \,  \, f(\omega_{L} \zeta_{L^c}) \right] \\
                    &= \int \, G_{\Lambda, z, \beta}(\mathrm{d}\zeta) \, g_{\Lambda, \Delta\setminus L}(\zeta) \left[\int \, G_{L, z, \beta, \zeta}(\mathrm{d}\omega) \,  \, f(\omega_{L} \zeta_{L^c}) \right] ,
                \end{align*} 
                as $G_{L, z, \beta, \zeta}$ is the appropriate conditional distribution and $g_{\Lambda, \Delta\setminus L}$ only depends on the particles in $\Delta \setminus L$.
                
                \item[4]   
                Using the fact that $0 \leq f \leq 1$ and that \(f\) just depends on the particles inside \(\Lambda_k\),
                \begin{align*}
                    &\bigg\vert \int \, G_{\Lambda, z, \beta}(\mathrm{d}\zeta) \, g_{\Lambda, \Delta\setminus L}(\zeta) \, \int \, G_{L, z, \beta, \zeta}(\mathrm{d}\omega) \, f(\omega_{L} \zeta_{L^c}) \\
                    &\qquad
                    - \int \, G_{\Lambda, z, \beta}(\mathrm{d}\zeta) \, g_{\Lambda, \Delta\setminus L}(\zeta) \, \int \, G_{L, z, \beta, \zeta_{\Lambda_k}}(\mathrm{d}\omega) \, f(\omega_{L} \zeta_{\Lambda_{k}\setminus L}) \bigg\vert \\
                    &\leq \int \, G_{\Lambda, z, \beta}(\mathrm{d}\zeta) \, g_{\Lambda, \Delta\setminus L}(\zeta) \left\vert \int \, G_{L, z, \beta, \zeta}(\mathrm{d}\omega) \, f(\omega_{L} \zeta_{L^c}) - \int \, G_{L, z, \beta, \zeta_{\Lambda_k}}(\mathrm{d}\omega) \, f(\omega_{L} \zeta_{\Lambda_{k}\setminus L}) \right\vert \\
                    &\leq \int \, G_{\Lambda, z, \beta}(\mathrm{d}\zeta) \, g_{\Lambda, \Delta\setminus L}(\zeta) \, \big\Vert G_{L, z, \beta, \zeta} - G_{L, z, \beta, \zeta_{\Lambda_k}} \big\Vert_{\text{TV}}.
                \end{align*}
                Now, using the fact that $G_{\Lambda\setminus(\Delta\setminus L), z, \beta, \zeta}$ is the conditional distribution of $G_{\Lambda, z, \beta}$ given the configuration $\zeta_{\Delta\setminus L}$ in $\Delta\setminus L$, as well as that $g_{\Lambda, \Delta\setminus L}$ just depends on the particles in $\Delta\setminus L$ and that $g_{\Lambda, \Delta\setminus L}$ is the density of $P$ with right to $G_{\Lambda, z, \beta}$ for events in $\Delta \setminus L$, we have that
                \begin{align*}
                    &\int \, G_{\Lambda, z, \beta}(\mathrm{d}\zeta) \, g_{\Lambda, \Delta\setminus L}(\zeta) \, \big\Vert G_{L, z, \beta, \zeta} - G_{L, z, \beta, \zeta_{\Lambda_k}} \big\Vert_{\text{TV}}\\
                    &= \int \, G_{\Lambda, z, \beta}(\mathrm{d}\zeta) \,\\
                    &\quad
                    \left[\int \, G_{\Lambda \setminus (\Delta \setminus L), z, \beta, \zeta}(\mathrm{d}\omega) \,  g_{\Lambda, \Delta\setminus L}(\zeta_{\Delta \setminus L}) \, \big\Vert G_{L, z, \beta, \omega_{\Lambda\setminus(\Delta\setminus C)}\zeta_{\Delta \setminus L}} - G_{L, z, \beta, (\omega_{\Lambda\setminus(\Delta\setminus C)}\zeta_{\Delta \setminus L})_{\Lambda_k}} \big\Vert_{\text{TV}} \right]\\
                    &= \int \, G_{\Lambda, z, \beta}(\mathrm{d}\zeta) \,\\
                    &\quad
                    g_{\Lambda, \Delta\setminus L}(\zeta) \left[\int \, G_{\Lambda \setminus (\Delta \setminus L), z, \beta, \zeta}(\mathrm{d}\omega)  \, \big\Vert G_{L, z, \beta, \omega_{\Lambda\setminus(\Delta\setminus C)}\zeta_{\Delta \setminus L}} - G_{L, z, \beta, (\omega_{\Lambda\setminus(\Delta\setminus C)}\zeta_{\Delta \setminus L})_{\Lambda_k}} \big\Vert_{\text{TV}} \right]\\
                    &= \int \, P(\mathrm{d}\zeta) \, \left[\int \, G_{\Lambda \setminus (\Delta \setminus L), z, \beta, \zeta}(\mathrm{d}\omega)  \, \big\Vert G_{L, z, \beta, \omega_{\Lambda\setminus(\Delta\setminus C)}\zeta_{\Delta \setminus L}} - G_{L, z, \beta, (\omega_{\Lambda\setminus(\Delta\setminus C)}\zeta_{\Delta \setminus L})_{\Lambda_k}} \big\Vert_{\text{TV}} \right]\\
                    &= \int \, PG_{\Lambda \setminus (\Delta \setminus L)}(\mathrm{d}\zeta) \, \big\Vert G_{L, z, \beta, \zeta} - G_{L, z, \beta, \zeta_{\Lambda_k}} \big\Vert_{\text{TV}}.
                \end{align*}
                By \Cref{le_total_variation_distance}, $\int \, PG_{\Lambda \setminus (\Delta \setminus L)}(\mathrm{d}\zeta) \, \big\Vert G_{L, z, \beta, \zeta} - G_{L, z, \beta, \zeta_{\Lambda_k}} \big\Vert_{\text{TV}} < \delta$.
                Therefore,
                \begin{align*}
                &\int \, G_{\Lambda, z, \beta}(\mathrm{d}\zeta) \, g_{\Lambda, \Delta\setminus L}(\zeta) \, \left[\int \, G_{L, z, \beta, \zeta}(\mathrm{d}\omega) \, f(\omega_{L} \zeta_{L^c}) \right] \\
                &\approx \int \, G_{\Lambda, z, \beta}(\mathrm{d}\zeta) \, g_{\Lambda, \Delta\setminus L}(\zeta) \, \left[\int \, G_{L, z, \beta, \zeta_{\Lambda_k}}(\mathrm{d}\omega) \, f(\omega_{L} \zeta_{\Lambda_{k}\setminus L}) \right].
                \end{align*}

                \item[5] Since $g_{\Lambda, \Delta\setminus L}(\zeta) = \frac{\mathrm{d} P_{\Delta\setminus L}(\zeta)}{\mathrm{d} G_{\Lambda, \Delta\setminus L, z, \beta}(\zeta)}$, we have 
                \begin{align*}
                    &\int \, G_{\Lambda, z, \beta}(\mathrm{d}\zeta) \, g_{\Lambda, \Delta\setminus L}(\zeta) \, \left[\int \, G_{L, z, \beta, \zeta_{\Lambda_k}}(\mathrm{d}\omega) \, f(\omega_{L} \zeta_{\Lambda_{k}\setminus L}) \right]\\
                    &= \int \, P(\mathrm{d}\zeta) \, \left[\int \, G_{L, z, \beta, \zeta_{\Lambda_k}}(\mathrm{d}\omega) \, f(\omega_{L} \zeta_{\Lambda_{k}\setminus L}) \right],
                \end{align*}
                as the corresponding integrand only depends on particles in $\Lambda_k\setminus L \subset \Delta\setminus L$.

                \item[6] This is true since, using the fact that $0 \leq f \leq 1$, we obtain
                \begin{align*}
                    &\left\vert \int \, P(\mathrm{d}\zeta) \, \left[\int \, G_{L, z, \beta, \zeta_{\Lambda_k}}(\mathrm{d}\omega) \, f(\omega_{L} \zeta_{\Lambda_{k}\setminus L}) \right] - \int \, P(\mathrm{d}\zeta) \, \left[\int \, G_{L, z, \beta, \zeta}(\mathrm{d}\omega) \, f(\omega_{L} \zeta_{\Lambda_{k}\setminus L}) \right] \right\vert \\
                    &\leq \int \, P(\mathrm{d}\zeta) \, \left\vert \int \, G_{L, z, \beta, \zeta_{\Lambda_k}}(\mathrm{d}\omega) \, f(\omega_{L} \zeta_{\Lambda_{k}\setminus L}) - \int \, G_{L, z, \beta, \zeta}(\mathrm{d}\omega) \, f(\omega_{L} \zeta_{\Lambda_{k}\setminus L})  \right\vert \\
                    &\leq \int \, P(\mathrm{d}\zeta) \, \left\Vert G_{L, z, \beta, \zeta_{\Lambda_k}} - G_{L, z, \beta, \zeta} \right\Vert_{\text{TV}}
                    < \delta,
                \end{align*} where the last inequality is again an application of \Cref{le_total_variation_distance}.
            \end{description}
            Putting the steps $\mathbf{(1)-(6)}$ together yields \eqref{eq:DLR_estimate}, as desired.
        \end{proof}

\subsection*{Acknowledgments}
BJ and JK gratefully received support by the Leibniz Association within the Leibniz Junior Research Group on \textit{Probabilistic Methods for Dynamic Communication Networks} as part of the Leibniz Competition (grant no.\ J105/2020).
AZ and BJ gratefully received support from Deutsche Forschungsgemeinschaft (DFG) through DFG Project no.\ 422743078 in the context of SPP 2265.

\bibliographystyle{alpha}    
\bibliography{references}

\end{document}